\documentclass[12pt,a4paper]{article}
\usepackage{amsfonts,amsmath,mathrsfs,amssymb, indentfirst}
\usepackage{rotating}
\usepackage{array}
\usepackage{tikz}
\newtheorem{theorem}{Theorem}[section]
\newtheorem{lemma}[theorem]{Lemma}

\newtheorem{defin}{Definition}[section]
\newtheorem{examp}{Example}[section]
\newtheorem{remark}{Remark}[section]
\newtheorem{cor}[theorem]{Corollary}

\begin{document}

\title {Block-diagonalized rigidity matrices of symmetric frameworks and applications}

\author{Bernd Schulze\footnote{Preparation of this manuscript was supported, in part, under a grant from NSERC (Canada), and final preparation occured at the TU Berlin with support of the DFG Research Unit 565 `Polyhedral Surfaces'.}\\ Inst. Mathematics, MA 6-2\\ TU Berlin\\ D-10623 Berlin, Germany
}

\maketitle

\begin{abstract}
In this paper, we give a complete self-contained proof that the rigidity matrix of a symmetric bar and joint framework (as well as its transpose) can be transformed into a block-diagonalized form using techniques from group representation theory. This theorem is basic to a number of useful and interesting results concerning the rigidity and flexibility of symmetric frameworks. As an example, we use this theorem to prove a generalization of the symmetry-extended version of Maxwell's rule given in \cite{FGsymmax} which can be applied to both injective and non-injective realizations in all dimensions.
\end{abstract}

 \section{Introduction}

It is a common method in engineering, physics, and chemistry to apply techniques from group representation theory  to the analysis of symmetric structures (see, for example, \cite{FG1, FG2, FG3, FG4, KG2, KG3}). In particular, some recent papers have used these techniques to gain insight into the rigidity properties of symmetric frameworks consisting of rigid bars and flexible joints \cite{cfgsw, FGsymmax, KG1, KG2}.\\\indent One of the fundamental observations resulting from this approach for studying the rigidity of symmetric bar and joint frameworks is due to R.D. Kangwai and S.D. Guest (\cite{KG2}):
given a symmetric framework $(G,p)$ and a non-trivial subgroup $S$ of its point group, there are techniques to block-diagonalize the rigidity matrix of $(G,p)$ (as well as its transpose) into submatrix blocks in such a way that each block corresponds to an irreducible representation of $S$. Using such a block-diagonalization, the (first-order) rigidity analysis of $(G,p)$ can be broken up into `symmetric' subproblems, where each subproblem considers the relationship between external forces on the joints of $(G,p)$ and resulting internal forces in the bars of $(G,p)$ that share certain symmetry properties. A number of interesting and useful results concerning the rigidity of symmetric frameworks are based on this method.\\\indent
However, since the main focus of the work in \cite{KG2}, as well as in subsequent papers such as \cite{cfgsw}, \cite{FGsymmax}, or \cite{KG1}, lies on applications in engineering and chemistry, many of these results are not presented with a mathematically precise formulation nor with a complete mathematical foundation.\\\indent
In this paper, we establish two major results. First, in Section 3, we use the mathematical foundation we established in \cite{BS1} to give a complete proof for the fact that the rigidity matrix of a symmetric framework (as well as its transpose) can be block-diagonalized in the way described above. Fundamental to this proof are our mathematically explicit definitions for the `external' and `internal' representation which were introduced in \cite{FGsymmax} and \cite{KG2} only by means of an example, and Lemma \ref{replemma} which establishes the key connection between these two representations.\\\indent
Secondly, in Section 4, we apply the results of Section 3 to give a detailed mathematical proof for the symmetry-extended version of Maxwell's rule given in \cite{FGsymmax}. This rule provides further necessary conditions (in addition to Maxwell's original condition from 1864 \cite{bibmaxwell}) for a symmetric framework to be isostatic (i.e., minimal infinitesimally rigid). While the symmetry-extended version of Maxwell's rule, as formulated in \cite{FGsymmax}, is only applicable to $2$- or $3$-dimensional frameworks with injective configurations (see \cite{BS1} for details), we establish a more general result in this paper, namely a rule that can be applied to both injective and non-injective realizations in \emph{all} dimensions. The proof of this result is based on Theorem \ref{Max} which in turn relies on the fact that the rigidity matrix of a symmetric framework can be block-diagonalized as described in Section 3. \\\indent An alternate approach to proving the symmetry-extended version of Maxwell's rule in \cite{FGsymmax}, as well as various generalizations of this rule to other types of geometric constraint systems, is given by J.C. Owen and
S.C. Power in \cite{owen}. \\\indent In order to apply the symmetry-extended version of Maxwell's rule to a given framework $(G,p)$, it is necessary to determine the dimensions of the subspaces of infinitesimal rigid motions of $(G,p)$ that are invariant under the external representation. While in \cite{FGsymmax}, the question of how to find the dimensions of these subspaces is only briefly addressed and not answered completely from a mathematical point of view (in particular, for all frameworks in dimensions higher than 3, this question is not addressed at all), in Section 4 of this paper, we explain in detail how to determine the dimensions of these subspaces for an arbitrary-dimensional framework.\\\indent
Since in \cite{FGsymmax} and \cite{KG2}, the rigidity properties of a symmetric framework are studied from both the kinematic and static point of view simultaneously, we develop the corresponding mathematical theory in this paper in the same manner.\\\indent
In \cite{cfgsw}, the symmetry-extended version of Maxwell's rule is used to show that a symmetric isostatic framework in 2D or 3D must obey some very simply stated restrictions on the number of structural elements that are `fixed' by various symmetry operations of the framework. Since the work in \cite{cfgsw} is based entirely on the symmetry-extended version of Maxwell's rule, the results of the present paper implicitly provide proofs for these results (and analogous extensions of these results to higher dimensions) as well.\\\indent
As shown in \cite{BS4, BS3, BS2}, symmetrized versions of some other famous theorems in rigidity theory, such as Laman's Theorem (\cite{graver, gss}) or the theorem on the equivalence of finite and infinitesimal rigidity for generic realizations (\cite{asiroth}), can also be established using and extending the results of this paper.

\section{Definitions and preliminaries}
\subsection{Introduction to infinitesimal and static rigidity}

\subsubsection{Infinitesimal rigidity}

We begin with a brief introduction to infinitesimal rigidity of bar and joint frameworks.

All graphs considered in this paper are finite graphs without loops or multiple edges. The \emph{vertex set} of a graph $G$ is denoted by $V(G)$ and the \emph{edge set} of $G$ is denoted by $E(G)$. Two vertices $u \ne v$ of $G$ are said to be \emph{adjacent} if $\{u,v\}\in E(G)$, and \emph{independent} otherwise.

\begin{defin}
\label{framework}
\emph{\cite{gss, W1, W2} A \emph{framework} (in $\mathbb{R}^{d}$) is a pair $(G,p)$, where $G$ is a graph and $p: V(G)\to \mathbb{R}^{d}$ is a map with the property that $p(u) \neq p(v)$ for all $\{u,v\} \in E(G)$. We also say that $(G,p)$
is a $d$-dimensional \emph{realization} of the \emph{underlying graph} $G$.}
\end{defin}

Given the vertex set $V(G)=\{v_{1},\ldots,v_{n}\}$ of a graph $G$ and a map $p:V(G)\to \mathbb{R}^{d}$, it is often useful to identify $p$ with a vector in $\mathbb{R}^{dn}$ by using the order on $V(G)$. In this case we also refer to $p$ as a \emph{configuration} of $n$ points in $\mathbb{R}^{d}$.

A \emph{joint} of a $d$-dimensional framework $(G,p)$ is an ordered pair $\big(v,p(v)\big)$, where $v \in V(G)$, and a \emph{bar} of $(G,p)$ is an unordered pair $\big\{\big(u,p(u)\big),\big(v,p(v)\big)\big\}$ of joints of $(G,p)$, where $\{u,v\} \in E(G)$. We define $\|p(u)-p(v)\|$ to be the \emph{length} of the bar $\big\{\big(u,p(u)\big),\big(v,p(v)\big)\big\}$, where $\|p(u)-p(v)\|$ is defined by the canonical inner product on $\mathbb{R}^{d}$.\\\indent
Note that the map $p$ of $(G,p)$ can possibly be non-injective, that is, two distinct joints $\big(u,p(u)\big)$ and $\big(v,p(v)\big)$ of $(G,p)$ may be located at the same point $p(u)=p(v)$ in $\mathbb{R}^{d}$, provided that $u$ and $v$ are independent vertices of $G$. However, if $\{u,v\} \in E(G)$, then $p(u)\neq p(v)$, so that every bar $\big\{\big(u,p(u)\big),\big(v,p(v)\big)\big\}$ of $(G,p)$ has a strictly positive length.

\begin{defin}
\label{infinmotion}
\emph{Let $(G,p)$ be a framework in $\mathbb{R}^{d}$ with $V(G)=\{v_{1},v_{2},\ldots, v_{n}\}$. An \emph{infinitesimal motion} of $(G,p)$ is a function $u: V(G)\to \mathbb{R}^{d}$ such that
\begin{equation}
\label{infinmotioneq}
\big(p(v_{i})-p(v_{j})\big)\cdot \big(u(v_{i})-u(v_{j})\big)=0 \quad\textrm{ for all } \{v_{i},v_{j}\} \in E(G)\textrm{.}\end{equation}
An infinitesimal motion $u$ of $(G,p)$ is an \emph{infinitesimal rigid motion} if there exists a family of differentiable functions $P_{i}:[0,1]\to \mathbb{R}^{d}, \, i=1,2,\ldots, n$, with $P_{i}(0)=p(v_{i})$ for all $i$ and $\|P_{i}(t)-P_{j}(t)\|=\|p(v_{i})-p(v_{j})\|$ for all $t\in [0,1]$ and all $1\leq i < j\leq n$, such that $u(v_{i})=P'_{i}(0)$ for all $i$. Otherwise $u$ is called an \emph{infinitesimal flex} of $(G,p)$.\\\indent
$(G,p)$ is said to be \emph{infinitesimally rigid} if every infinitesimal motion of $(G,p)$ is an infinitesimal rigid motion. Otherwise $(G,p)$ is said to be \emph{infinitesimally flexible}. See \cite{gss, W1}, for example, for more details.}
\end{defin}

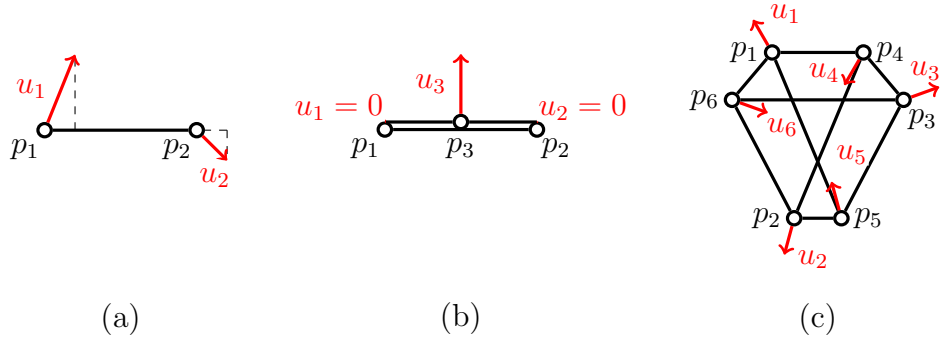
\begin{figure}[htp]
\begin{center}
\begin{tikzpicture}[very thick,scale=1]
\tikzstyle{every node}=[circle, draw=black, fill=white, inner sep=0pt, minimum width=5pt];
        \path (0,0) node (p1) [label = below left: $p_{1}$] {} ;
        \path (2,0) node (p2) [label = below left: $p_2$] {} ;
        \draw (p1)  --  (p2);
        \draw [dashed, thin] (0.4,0) -- (0.4,1);
        \draw [dashed, thin] (2.4,0) -- (2.4,-0.4);
        \draw [dashed, thin] (2.4,0) -- (p2);
        \draw [->, red] (p1) -- node [draw=white, left=4pt] {$u_{1}$} (0.4,1);
        \draw [->, red] (p2) -- node [draw=white, below =4pt] {$u_2$} (2.4,-0.4);
        \node [draw=white, fill=white] (a) at (1,-2.5) {(a)};
        \end{tikzpicture}
        \hspace{0.5cm}
        \begin{tikzpicture}[very thick,scale=1]
\tikzstyle{every node}=[circle, draw=black, fill=white, inner sep=0pt, minimum width=5pt];
    \path (0,0) node (p1) [label = below left: $p_1$] {} ;
    \path (2,0) node (p2) [label = below right: $p_2$] {} ;
    \node (p3) at (1,0.1) {};
    \node [draw=white, fill=white] (labelp3) at (1,-0.25) {$p_{3}$};
    \draw (p1) -- (p2);
    \draw (0,0.1) -- (p3);
    \draw (p3) -- (2,0.1);
    \draw [->, red] (p3) -- node [draw=white, left=4pt] {$u_{3}$} (1,1);
    \draw [ red] (p1) -- node [rectangle, draw=white, above left] {$u_{1}=0$} (p1);
    \draw [ red] (p2) -- node [rectangle, draw=white, above right] {$u_{2}=0$} (p2);
    \node [draw=white, fill=white] (b) at (1,-2.5) {(b)};
    \end{tikzpicture}
    \hspace{0.5cm}
        \begin{tikzpicture}[very thick,scale=1]
\tikzstyle{every node}=[circle, draw=black, fill=white, inner sep=0pt, minimum width=5pt];
    \path (160:1.2cm) node (p6) [label = left: $p_6$] {} ;
    \path (120:1.2cm) node (p1) [label = left: $p_1$] {} ;
    \path (255:1.2cm) node (p2) [label = left: $p_2$] {} ;
     \path (20:1.2cm) node (p3) [label = below right: $p_3$] {} ;
    \path (60:1.2cm) node (p4) [label = right: $p_4$] {} ;
    \path (285:1.2cm) node (p5) [label = right: $p_5$] {} ;
     \draw (p1) -- (p4);
     \draw [->, red] (p6) --  (160:0.7)node[draw=white, below right=0.3pt]{$u_{6}$};
     \draw (p1) -- (p5);
     \draw (p1) -- (p6);
     \draw (p2) -- (p4);
     \draw (p2) -- (p5);
     \draw (p2) -- (p6);
     \draw (p3) -- (p4);
     \draw (p3) -- (p5);
     \draw (p3) -- (p6);
     \draw [->, red] (p1) -- node[rectangle, draw=white, above right=3pt] {$u_{1}$}(120:1.7);
     \draw [->, red] (p2) -- node[rectangle, draw=white, below right=3pt] {$u_{2}$}(255:1.7);
     \draw [->, red] (p3) -- node[rectangle, draw=white, above=4pt] {$u_{3}$}(20:1.7);
     \draw [->, red] (p4)  -- node[rectangle, draw=white, left=4pt] {$u_{4}$} node[rectangle, draw=white, below=27pt] {$u_{5}$}(60:0.7);
     \draw [->, red] (p5) -- (285:0.7);
     \node [draw=white, fill=white] (c) at (0,-2.5) {(c)};
        \end{tikzpicture}
\end{center}
\caption{\emph{The arrows indicate the non-zero displacement vectors of an infinitesimal rigid motion \emph{(a)} and infinitesimal flexes \emph{(b, c)} of frameworks in $\mathbb{R}^2$.}}
\label{inmo}
\end{figure}

Note that an infinitesimal motion of a framework $(G,p)$ is a set of displacement vectors, one at each joint, that neither stretch nor compress the bars of $(G,p)$ at first order. More precisely, condition (\ref{infinmotioneq}) says that for every edge $\{v_{i},v_{j}\} \in E(G)$, the projections of $u(v_{i})$ and $u(v_{j})$ onto the line through $p(v_{i})$ and $p(v_{j})$ have the same direction and the same length (see also Figure \ref{inmo}) \cite{W1, W2}.

\begin{remark}
\label{infinmotionrem}
\emph{Let $G$ be a graph with $V(G)=\{v_{1},v_{2},\ldots, v_{n}\}$ and let $u$ be an infinitesimal motion of a $d$-dimensional realization $(G,p)$ of $G$. If $\big(p(v_{i})-p(v_{j})\big)\cdot \big(u(v_{i})-u(v_{j})\big)\neq 0$ for some $\{v_{i},v_{j}\} \notin E(G)$, then $u$ is an infinitesimal flex of $(G,p)$. If the points $p(v_{1}),\ldots ,p(v_{n})$ span all of $\mathbb{R}^d$ (in an affine sense), then the converse also holds, i.e., in this case, $u$ is an infinitesimal flex of $(G,p)$ if and only if $\big(p(v_{i})-p(v_{j})\big)\cdot \big(u(v_{i})-u(v_{j})\big)\neq 0$ for some $\{v_{i},v_{j}\} \notin E(G)$ or equivalently, $u$ is an infinitesimal rigid motion of $(G,p)$ if and only if $\big(p(v_{i})-p(v_{j})\big)\cdot \big(u(v_{i})-u(v_{j})\big)=0$ for all $1\leq i < j\leq n$.}
\end{remark}

From now on, when we say that a set of points spans a space, then this will always be in the affine sense.

For a framework $(G,p)$ whose underlying graph $G$ has a vertex set that is indexed from 1 to $n$, say $V(G)=\{v_{1},v_{2},\ldots ,v_{n}\}$, we will frequently denote $p(v_{i})$ by $p_{i}$ for $i=1,2,\ldots, n$. The $k^{th}$ component of a vector $x$ is denoted by $(x)_{k}$.

The equations in Definition \ref{infinmotion} form a system of linear equations whose corresponding matrix is called the rigidity matrix. This matrix is fundamental in the study of both infinitesimal and static rigidity \cite{graver, gss, W1, W2}.

\begin{defin}
\emph{ Let $G$ be a graph with $V(G)=\{v_{1},v_{2},\ldots,v_{n}\}$ and let $p:V(G)\to \mathbb{R}^{d}$. The
\emph{rigidity matrix} of $(G,p)$ is the $|E(G)| \times dn$ matrix  \\\indent \begin{displaymath} \mathbf{R}(G,p)=\left(
\begin{array} {ccccccccccc }
& & & & & \vdots & & & & & \\
0 & \ldots & 0 & p_{i}-p_{j}&0 &\ldots &0 & p_{j}-p_{i} &0 &\ldots &
0\\ & & & & & \vdots & & & & &\end{array}
\right)\textrm{,}\end{displaymath} that is, for each edge $\{v_{i},v_{j}\}\in E(G)$, $\mathbf{R}(G,p)$ has the row with
$(p_{i}-p_{j})_{1},\ldots,(p_{i}-p_{j})_{d}$ in the columns $d(i-1)+1,\ldots,di$, $(p_{j}-p_{i})_{1},\ldots,(p_{j}-p_{i})_{d}$ in
the columns $d(j-1)+1,\ldots,dj$, and $0$ elsewhere.}
\end{defin}

\begin{remark}
\emph{The rigidity matrix is defined for arbitrary pairs $(G,p)$, where $G$ is a graph and $p:V(G)\to \mathbb{R}^{d}$ is a map. If $(G,p)$ is not a framework, then there exists a pair of adjacent vertices of $G$ that are mapped to the same point in $\mathbb{R}^{d}$ under $p$ and every such edge of $G$ gives rise to a zero-row in $\mathbf{R}(G,p)$.}
\end{remark}

If we identify an infinitesimal motion of a $d$-dimensional framework $(G,p)$ with a column vector in $\mathbb{R}^{d|V(G)|}$ (by using the order on $V(G)$), then the kernel of $\mathbf{R}(G,p)$ is the space of infinitesimal motions of $(G,p)$. It is known that the infinitesimal rigid motions arising from $d$ translations and $\binom{d}{2}$ rotations of $\mathbb{R}^{d}$ form a basis of the space of infinitesimal rigid motions of $(G,p)$, provided that the points $p_{1},\ldots ,p_{n}$ span an affine subspace of $\mathbb{R}^{d}$ of dimension at least $d-1$ \cite{gss, W1}. Thus, for such a framework $(G,p)$, we have $\textrm{nullity }\big(\mathbf{R}(G,p)\big)\geq d+\binom{d}{2}=\binom{d+1}{2}$ and $(G,p)$ is infinitesimally rigid if and only if $\textrm{nullity } \big(\mathbf{R}(G,p)\big)=\binom{d+1}{2}$ or equivalently, $\textrm{rank }\big(\mathbf{R}(G,p)\big)=d |V(G)| - \binom{d+1}{2}$.

\begin{theorem}
\label{infinrigaff}\cite{asiroth}
A framework $(G,p)$ in $\mathbb{R}^d$ is infinitesimally rigid if and only if either $\textrm{rank }\big(\mathbf{R}(G,p)\big)=d |V(G)| - \binom{d+1}{2}$ or $G$ is a complete graph $K_n$ and the points $p(v)$, $v\in V(G)$, are affinely independent.
\end{theorem}

\subsubsection{Static rigidity}

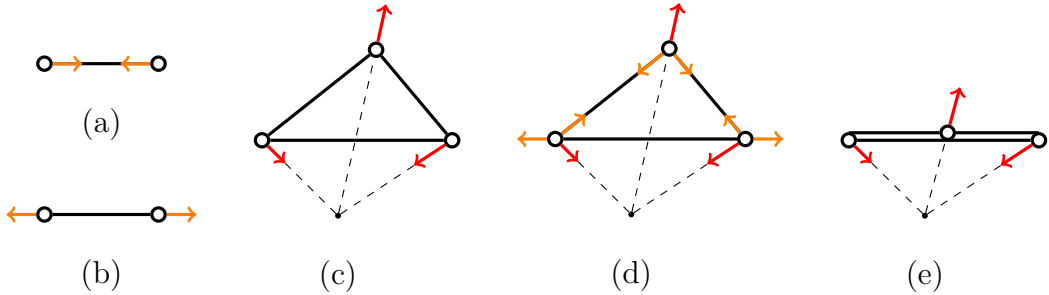
\begin{figure}[htp]
\begin{center}
\begin{tikzpicture}[very thick,scale=1]
\tikzstyle{every node}=[circle, draw=black, fill=white, inner sep=0pt, minimum width=5pt];
    \path (0:0cm) node (p1) {};
    \path (0:1.5cm) node (p2) {};
    \draw (p1) -- (p2);
    \draw [->, orange] (p1) -- (0.5,0);
    \draw [->, orange] (p2) -- (1,0);

    \node [draw=white, fill=white] (a) at (0.75,-0.8) {(a)};

    \path (270:2cm) node (p3) {};
    \node (p4) at (1.5,-2) {};
    \draw (p3) -- (p4);
    \draw [->, orange] (p3) -- (-0.5,-2);
    \draw [->, orange] (p4) -- (2,-2);

    \node [draw=white, fill=white] (b) at (0.75,-2.8) {(b)};
    \end{tikzpicture}
    \hfill
    \begin{tikzpicture}[very thick,scale=1]
    \tikzstyle{every node}=[circle, draw=black, fill=white, inner sep=0pt, minimum width=5pt];
    \node (p5) at (0,0) {};
    \node (p6) at (2.5,0) {};
    \node (p7) at (1.5,1.2) {};
    \draw (p5) -- (p6);
    \draw (p6) -- (p7);
    \draw (p7) -- (p5);
    \path (1,-1) node [circle, draw=black, fill=black, minimum width=0.5pt](o) {};
    \draw [dashed, thin] (o) -- (p5);
    \draw [dashed, thin] (o) -- (p6);
    \draw [dashed, thin] (o) -- (p7);
    \draw [->, red] (p5) -- (0.3,-0.3);
    \draw [->, red] (p6) -- (2,-0.33);
    \draw [->, red] (p7) -- (1.63,1.8);
    \node [draw=white, fill=white] (c) at (1,-1.8) {(c)};
    \end{tikzpicture}
    \hfill
    \begin{tikzpicture}[very thick,scale=1]
    \tikzstyle{every node}=[circle, draw=black, fill=white, inner sep=0pt, minimum width=5pt];
    \node (p5) at (0,0) {};
    \node (p6) at (2.5,0) {};
    \node (p7) at (1.5,1.2) {};
    \draw (p5) -- (p6);
    \draw (p6) -- (p7);
    \draw (p7) -- (p5);
    \path (1,-1) node [circle, draw=black, fill=black, minimum width=0.5pt](o) {};
    \draw [dashed, thin] (o) -- (p5);
    \draw [dashed, thin] (o) -- (p6);
    \draw [dashed, thin] (o) -- (p7);
    \draw [->, red] (p5) -- (0.3,-0.3);
    \draw [->, red] (p6) -- (2,-0.33);
    \draw [->, red] (p7) -- (1.63,1.8);
    \draw [->, orange] (p5) -- (-0.5,0);
    \draw [->, orange] (p5) -- (0.4,0.32);
    \draw [->, orange] (p6) -- (3,0);
    \draw [->, orange] (p6) -- (2.23,0.32);
    \draw [->, orange] (p7) -- (1.08,0.86);
    \draw [->, orange] (p7) -- (1.79,0.85);
    \node [draw=white, fill=white] (d) at (1,-1.8) {(d)};
    \end{tikzpicture}
    \hfill
    \begin{tikzpicture}[very thick,scale=1]
    \tikzstyle{every node}=[circle, draw=black, fill=white, inner sep=0pt, minimum width=5pt];
    \node (p5) at (0,0) {};
    \node (p6) at (2.5,0) {};
    \node (p7) at (1.3,0.1) {};
    \draw (p5) -- (p6);
    \draw (2.5,0.1) -- (p7);
    \draw (p7) -- (0,0.1);
    \path (1,-1) node [circle, draw=black, fill=black, minimum width=0.5pt](o) {};
    \draw [dashed, thin] (o) -- (p5);
    \draw [dashed, thin] (o) -- (p6);
    \draw [dashed, thin] (o) -- (p7);
    \draw [->, red] (p5) -- (0.3,-0.3);
    \draw [->, red] (p6) -- (2,-0.33);
    \draw [->, red] (p7) -- (1.46,0.7);
    \node [draw=white, fill=white] (e) at (1,-1.8) {(e)};
    \end{tikzpicture}
\end{center}
\caption{\emph{\emph{(a), (b)} The arrows indicate a tension \emph{(a)} and a compression \emph{(b)} in a bar. \emph{(c)} An equilibrium load on a non-degenerate triangle. This load can be resolved by the triangle as shown in \emph{(d)}. \emph{(e)} An unresolvable equilibrium load on a degenerate triangle: for any joint of this framework, tensions or compressions in the bars  cannot reach an equilibrium with the load vector at this joint.}}
\label{fig:eqload}
\end{figure}

We now also give a brief introduction to the static approach to rigidity. The intuitive test for static rigidity of a framework $(G,p)$ is to apply an external load to $(G,p)$ (i.e., a set of forces, one to each joint) and investigate whether there exists a set of tensions and compressions in the bars of $(G,p)$ that reach an equilibrium with this load at the joints (see also Figure \ref{fig:eqload}). Of course only loads which do not correspond to a translation or rotation of space can possibly be resolved in this way.

\begin{defin}
\label{load}
\emph{\cite{CW, taywhit, W4, W1} Let $(G,p)$ be a framework in $\mathbb{R}^{d}$ with $V(G)=\{v_{1},v_{2},\ldots,v_{n}\}$. A \emph{load} on $(G,p)$ is a function $l:V(G)\to \mathbb{R}^{d}$, where for $i=1,2,\ldots,n$, the vector $l(v_{i})$ represents a force applied to the joint $\big(v_{i},p_{i}\big)$ of $(G,p)$.\\\indent A load $l$ on $(G,p)$ is called an \emph{equilibrium load} if $l$ satisfies
\begin{itemize}
\item[$(i)$] $\sum_{i=1}^{n}l_{i}=0$;
\item[$(ii)$] $\sum_{i=1}^{n}\big((l_{i})_{j}(p_{i})_{k}- (l_{i})_{k}(p_{i})_{j}\big)=0$ \quad for all $1\leq j< k\leq d$,
\end{itemize}
where $l_{i}$ denotes the vector $l(v_{i})$ for each $i$.}
\end{defin}

The physical intuition for conditions $(i)$ and $(ii)$ in Definition \ref{load} is
the following: condition $(i)$ rules out loads that would produce a translation of $(G,p)$ and $(ii)$ says that there is no net rotational twist of $(G,p)$.

\begin{defin}
\label{resolution}
\emph{\cite{CW, taywhit, W4, W1} Let $l$ be an equilibrium load on a framework $(G,p)$ in $\mathbb{R}^{d}$ with $V(G)=\{v_{1},v_{2},\ldots,v_{n}\}$. A \emph{resolution} of $l$ by $(G,p)$ (also called a \emph{stress} of $(G,p)$) is a function $\omega:E(G)\to \mathbb{R}$ such that at each joint $\big(v_{i},p_{i}\big)$ of $(G,p)$ we have
\begin{displaymath}
\sum_{j \textrm{ with }\{v_{i},v_{j}\}\in E(G)}\omega_{ij}(p_{i}-p_{j})+ l_{i}=0 \textrm{,}
\end{displaymath}
where $\omega_{ij}$ denotes $\omega(\{v_{i},v_{j}\})$ for all $\{v_{i},v_{j}\}\in E(G)$.}
\end{defin}

The scalars $\omega_{ij}$ represent tensions ($\omega_{ij}<0$) and compressions ($\omega_{ij}>0$) in the bars of $(G,p)$, so that the bar forces reach an equilibrium with $l_{i}$ at each joint $\big(v_{i},p_{i}\big)$.

\begin{defin}
\label{statrigid}
\emph{\cite{CW, taywhit, W4, W1} A framework $(G,p)$ is \emph{statically rigid} if every equilibrium load on $(G,p)$ has a resolution by $(G,p)$.}
\end{defin}

Note that if we identify $l$ and $\omega$ with a column vector in $\mathbb{R}^{dn}$ and $\mathbb{R}^{|E(G)|}$, respectively, then (after changing the sign of $l$) the equations in Definition \ref{resolution} can be written in a compact form in terms of the rigidity matrix $\mathbf{R}(G,p)$ as \begin{displaymath}\mathbf{R}(G,p)^{T} \omega =l\textrm{.}\end{displaymath} Let $(v_{h},p_{h})$ and $(v_{k},p_{k})$ be two joints of $(G,p)$. Then it is easy to see that the column vector $F_{hk}$, where \begin{displaymath}(F_{hk})^{T}=(0,\ldots, 0, p_{h}-p_{k}, 0 ,\ldots ,0 ,p_{k}-p_{h},0 ,\ldots ,0)\textrm{,}\end{displaymath} is an equilibrium load on $(G,p)$. Further, if $\{v_{i},v_{j}\}\in E(G)$, then $(F_{ij})^{T}$ is the row vector of $\mathbf{R}(G,p)$ that corresponds to $\{v_{i},v_{j}\}$ and $F_{ij}$ is clearly resolved by the bar $\{(v_{i},p_{i}),(v_{j},p_{j})\}$ of $(G,p)$. Note that if $(G,p)$ is statically rigid, then $F_{hk}$ has a resolution by $(G,p)$ for \emph{every} pair $(v_{h},p_{h}),(v_{k},p_{k})$ of joints of $(G,p)$ (even if $\{v_{h},v_{k}\} \notin E(G)$).\\\indent If the points $p_{1},\ldots, p_{n}$ span all of $\mathbb{R}^d$, then the converse also holds, since in this case, the vectors $F_{hk}$, $1\leq h<k\leq n$, generate the entire space of equilibrium loads on $(G,p)$ (see \cite{W4}). This space is a subspace of $\mathbb{R}^{d n}$ of dimension $d n -\binom{d+1}{2}$ (defined by the equations in Definition \ref{load}).\\\indent Thus, if we want to test such a framework $(G,p)$ for static rigidity, we need to investigate whether the rows of $\mathbf{R}(G,p)$ generate a space of dimension $d n -\binom{d+1}{2}$, that is, the entire space of equilibrium loads on $(G,p)$. In other words, we need to investigate whether \begin{displaymath}\textrm{rank }\big(\mathbf{R}(G,p)^{T}\big)=d n -\binom{d+1}{2}\textrm{.}\end{displaymath}

So, the essential information for both infinitesimal and static rigidity of a framework $(G,p)$ is comprised by the rigidity matrix $\mathbf{R}(G,p)$. While in infinitesimal rigidity, we investigate the column space and column rank of $\mathbf{R}(G,p)$, in static rigidity, we investigate the row space and row rank of $\mathbf{R}(G,p)$. In the light of these remarks, the following fundamental facts do not come as a surprise.

\begin{theorem}
\label{flexandload}\cite{RW1}
The load $F_{hk}$ on a framework $(G,p)$ has no resolution by $(G,p)$ if and only if there exists an infinitesimal motion $u$ of $(G,p)$ with $(p_{h}-p_{k})\cdot (u_{h}-u_{k})\neq 0$.
\end{theorem}

\begin{theorem}
\label{statandinfequiv}\cite{henne, RW1}
A framework $(G,p)$ is infinitesimally rigid if and only if $(G,p)$ is statically rigid.
\end{theorem}

Theorem \ref{statandinfequiv} allows us to use the terms infinitesimally rigid and statically rigid interchangeably.

\subsection{Symmetry in frameworks}

We now introduce the necessary terms and definitions relating to symmetric frameworks.

First, recall that an \emph{automorphism} of a graph $G$ is a permutation $\alpha$ of $V(G)$ such that $\{u,v\}\in E(G)$ if and only if $\{\alpha(u),\alpha(v)\}\in E(G)$. The automorphisms of a graph $G$ form a group under composition which is denoted by $\textrm{Aut}(G)$.\\\indent Also, recall that an \emph{isometry} of $\mathbb{R}^{d}$ is a map $x:\mathbb{R}^{d}\to \mathbb{R}^{d}$ such that $\|x(a)-x(b)\|=\|a-b\|$ for all $a,b\in \mathbb{R}^{d}$.

\begin{defin}
\label{symop}
\emph{\cite{Hall, BS1} Let $(G,p)$ be a framework in $\mathbb{R}^{d}$.
A \emph{symmetry operation} of $(G,p)$ is an isometry $x$ of $\mathbb{R}^{d}$ such that for some $\alpha\in \textrm{Aut}(G)$, we have
$x\big(p(v)\big)=p\big(\alpha(v)\big)$ for all $v\in V(G)$.}
\end{defin}

The set of all symmetry operations of a given framework forms a group under composition. We adopt the following vocabulary from chemistry and crystallography:

\begin{defin} \emph{Let $(G,p)$ be a framework. Then the group which consists of all symmetry operations of $(G,p)$ is called the \emph{point group} of $(G,p)$.}
\end{defin}

It is well known that if $P$ is the point group of a $d$-dimensional framework $(G,p)$, then there exists a point $O$ in $\mathbb{R}^{d}$ which is fixed by every symmetry operation in $P$ \cite{BS1}. Since a translation of $(G,p)$ does not change the rigidity properties of $(G,p)$, we may assume wlog that this point $O$ is the origin of $\mathbb{R}^{d}$. It then follows that every symmetry operation of $(G,p)$ is an orthogonal linear transformation of $\mathbb{R}^{d}$, so that $P$ is a \emph{symmetry group}, i.e., a subgroup of the orthogonal group $O(\mathbb{R}^{d})$.\\\indent In this paper, the point group of every framework is assumed to be a symmetry group.

For the symmetry operations and symmetry groups of the $2$-dimensional frameworks given in the examples of this paper, we use the Schoenflies notation since it is one of the standard notations in the literature about symmetric structures. The three kinds of possible symmetry operations in dimension $2$ are the identity $Id$, rotations $C_{m}$ about the origin by an angle of $\frac{2\pi}{m}$, where $m\geq 2$, and reflections $s$ in lines through the origin. In the Schoenflies notation, this gives rise to the following families of possible symmetry groups in dimension 2: $\mathcal{C}_{1}$, $\mathcal{C}_{s}$, $\mathcal{C}_{m}$ and $\mathcal{C}_{mv}$, where $m\geq 2$. $\mathcal{C}_{1}$ denotes the trivial group which only contains the identity $Id$. $\mathcal{C}_{s}$ denotes any symmetry group in dimension 2 that consists of the identity $Id$ and a single reflection $s$. For $m\geq 2$, $\mathcal{C}_{m}$ denotes any cyclic symmetry group of order $m$ which is generated by a rotation $C_{m}$, and $\mathcal{C}_{mv} $ denotes any symmetry group in dimension 2 that is generated by a pair $\{C_{m},s\}$.

In order to symmetrize results in rigidity theory, we need an appropriate classification of symmetric frameworks. We use the following terminology (see also \cite{BS1, BS4, BS3, BS2}).

\begin{defin}
\label{symclass}
\emph{Let $G$ be a graph and $S$ be a symmetry group in dimension $d$. Then \emph{$\mathscr{R}_{(G,S)}$} is the set of all $d$-dimensional realizations of $G$ whose point group is either equal to $S$ or contains $S$ as a subgroup. An element of $\mathscr{R}_{(G,S)}$ is said to be a \emph{realization of the pair $(G,S)$}.}
\end{defin}

It follows directly from these definitions that if $(G,p)$ is a $d$-dimensional realization of a graph $G$ and $S$ is a symmetry group in dimension $d$, then $(G,p)\in \mathscr{R}_{(G,S)}$ if and only if there exists a map $\Phi:S\to \textrm{Aut}(G)$ such that
\begin{equation}\label{class} x\big(p(v)\big)=p\big(\Phi(x)(v)\big)\textrm{ for all } v\in V(G)\textrm{ and all } x\in S\textrm{.}\end{equation}

If a framework $(G,p)\in \mathscr{R}_{(G,S)}$ satisfies the equations in (\ref{class}) for the map $\Phi:S\to \textrm{Aut}(G)$, then
$(G,p)$ is said to be \emph{of type $\Phi$}. The set of all realizations of $(G,S)$ which are of type $\Phi$ is denoted by $\mathscr{R}_{(G,S,\Phi)}$.

Given a graph $G$ and a symmetry group $S$ in dimension $d$, different choices of types $\Phi:S\to \textrm{Aut}(G)$ frequently lead to very different geometric types of realizations of $(G,S)$. This is because a type $\Phi$ forces the joints and bars of a framework in $\mathscr{R}_{(G,S,\Phi)}$ to assume certain geometric positions in $\mathbb{R}^d$, as the following example demonstrates.

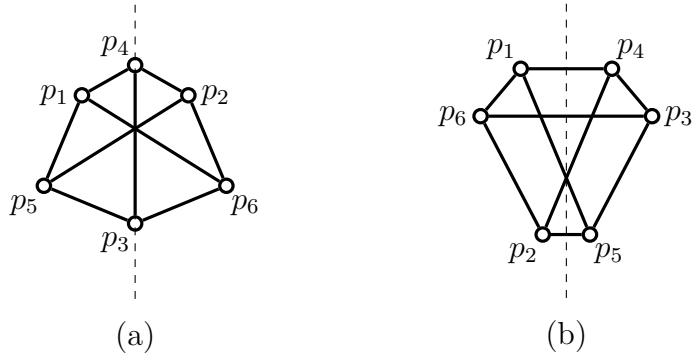
\begin{figure}[htp]
\begin{center}
\begin{tikzpicture}[very thick,scale=1]
\tikzstyle{every node}=[circle, draw=black, fill=white, inner sep=0pt, minimum width=5pt];
    \path (0.3,-0.5) node (p5) [label = below left: $p_{5}$] {} ;
    \path (1.5,-1) node (p3) [label = below left: $p_{3}$] {} ;
    \path (2.7,-0.5) node (p6) [label = below right: $p_{6}$] {} ;
   \path (0.8,0.7) node (p1) [label = left: $p_{1}$] {} ;
   \path (2.2,0.7) node (p2) [label = right: $p_{2}$] {} ;
   \path (1.5,1.1) node (p4) [label = above left: $p_{4}$] {} ;
   \draw (p1) -- (p4);
     \draw (p1) -- (p5);
     \draw (p1) -- (p6);
     \draw (p2) -- (p4);
     \draw (p2) -- (p5);
     \draw (p2) -- (p6);
     \draw (p3) -- (p4);
     \draw (p3) -- (p5);
     \draw (p3) -- (p6);
   \draw [dashed, thin] (1.5,-2) -- (1.5,2);
      \node [draw=white, fill=white] (a) at (1.5,-2.5) {(a)};
    \end{tikzpicture}
    \hspace{2cm}
        \begin{tikzpicture}[very thick,scale=1]
\tikzstyle{every node}=[circle, draw=black, fill=white, inner sep=0pt, minimum width=5pt];
    \path (160:1.2cm) node (p6) [label = left: $p_6$] {} ;
    \path (120:1.2cm) node (p1) [label = above left: $p_1$] {} ;
    \path (255:1.2cm) node (p2) [label = below left: $p_2$] {} ;
     \path (20:1.2cm) node (p3) [label = right: $p_3$] {} ;
    \path (60:1.2cm) node (p4) [label = above right: $p_4$] {} ;
    \path (285:1.2cm) node (p5) [label = below right: $p_5$] {} ;
     \draw (p1) -- (p4);
     \draw (p1) -- (p5);
     \draw (p1) -- (p6);
     \draw (p2) -- (p4);
     \draw (p2) -- (p5);
     \draw (p2) -- (p6);
     \draw (p3) -- (p4);
     \draw (p3) -- (p5);
     \draw (p3) -- (p6);
     \draw [dashed, thin] (0,-2) -- (0,2);
      \node [draw=white, fill=white] (b) at (0,-2.5) {(b)};
        \end{tikzpicture}
\end{center}
\caption{\emph{$2$-dimensional realizations of $(K_{3,3},\mathcal{C}_{s})$ of different types.}}
\label{K33types}
\end{figure}

\begin{examp}
\label{K33ex}
\emph{Figure \ref{K33types} shows two realizations of $(K_{3,3},\mathcal{C}_{s})$ of different types, where $K_{3,3}$ is the complete bipartite graph with partite sets $\{v_{1},v_{2},v_{3}\}$ and $\{v_{4},v_{5},v_{6}\}$, and $\mathcal{C}_{s}=\{Id,s\}$ is a symmetry group in dimension 2 generated by a reflection $s$. The framework in Figure \ref{K33types} (a) is a realization of $(K_{3,3},\mathcal{C}_{s})$ of type
$\Phi_{a}$, where $\Phi_{a}: \mathcal{C}_{s} \to \textrm{Aut}(K_{3,3})$ is defined by
\begin{eqnarray} \Phi_{a}(Id)& =& id\nonumber\\\Phi_{a}(s)&=&
(v_{1}\,v_{2})(v_{5}\,v_{6})(v_{3})(v_{4})\textrm{,}\nonumber
\end{eqnarray}
and the framework in Figure \ref{K33types} (b)
is a realization of $(K_{3,3},\mathcal{C}_{s})$ of type
$\Phi_{b}$, where $\Phi_{b}: \mathcal{C}_{s} \to \textrm{Aut}(K_{3,3})$ is defined by
\begin{eqnarray}\Phi_{b}(Id)& = &id\nonumber\\ \Phi_{b}(s) &=& (v_{1}\,v_{4})(v_{2}\,v_{5})(v_{3}\,v_{6})\textrm{.}\nonumber
\end{eqnarray}}
\end{examp}

As shown in \cite{BS1}, `almost all' realizations in a set of the form $\mathscr{R}_{(G,S,\Phi)}$ share the same infinitesimal rigidity properties. \\\indent For example, `almost all' realizations in $\mathscr{R}_{(K_{3,3},\mathcal{C}_{s},\Phi_{a})}$ are infinitesimally rigid, whereas all realizations in $\mathscr{R}_{(K_{3,3},\mathcal{C}_{s},\Phi_{b})}$ are infinitesimally flexible \cite {BS1, WW, W3}.\\\indent It is also shown in \cite{BS1} that if $(G,p)\in\mathscr{R}_{(G,S)}$ is an injective realization of $G$, then $(G,p)$ is of a unique type $\Phi$ and $\Phi$ must be a homomorphism. However, if $(G,p)$ is a non-injective realization, then $(G,p)$ may be of several types and a given type may not be a homomorphism. For details, we refer the reader to \cite{BS1, BS4}.\\\indent As we will see in this paper, the type $\Phi$ plays a key role in applying techniques from group representation theory to the analysis of a symmetric framework in a set of the form $\mathscr{R}_{(G,S)}$, whenever $\Phi$ is a homomorphism.

\begin{remark}\emph{Let $G$ be a graph, $S$ be a symmetry group in dimension $d$, and $\Phi$ be a map from $S$ to $\textrm{Aut}(G)$. Then it is easy to see that for any given $x\in S$, the set of all configurations $p$ of $n$ points in $\mathbb{R}^{d}$ that satisfy the equations in (\ref{class}) corresponding to $x$ is a linear subspace of $\mathbb{R}^{dn}$ (see also \cite{BS1}). We denote this subspace by $L_{x,\Phi}$. It follows that $U= \bigcap_{x \in S} L_{x,\Phi}$ is also a subspace of $\mathbb{R}^{dn}$. Note that $U$ is the space of all those (possibly non-injective) configurations $p$ of $n$ points in $\mathbb{R}^{d}$ with the property that either $(G,p)$ is a framework in $\mathscr{R}_{(G,S,\Phi)}$, or $p$ possesses the symmetry imposed by $S$ and $\Phi$, but there exists at least one edge $\{v_i,v_j\}$ in $E(G)$ with $p_i=p_j$.}
\end{remark}

\subsection{Basic definitions in group representation theory}

We need the following notions from group representation theory.

\begin{defin}
\emph{Let $S$ be a group and $V$ be an $n$-dimensional vector space over the field $F$. A \emph{linear representation} of $S$ with \emph{representation space} $V$ is a group homomorphism $H$ from $S$ to $GL(V)$, where $GL(V)$ denotes the group of all automorphisms of $V$. The dimension $n$ of $V$ is called the \emph{degree} of $H$.\\\indent Two linear representations $H_{1}:S\to GL(V_{1})$ and $H_{2}:S\to GL(V_{2})$ are said to be \emph{equivalent} if there exists an isomorphism $h:V_{1}\to V_{2}$ such that $h\circ H_{1}(x)\circ h^{-1}=H_{2}(x)$ for all $x\in S$.}
\end{defin}

\begin{defin}\emph{Let $S$ be a group, $V$ be a vector space over the field $F$ and $H:S\to GL(V)$ be a linear representation of $S$. A subspace $U$ of $V$ is said to be \emph{$H$-invariant} (or simply \emph{invariant} if $H$ is clear from the context) if $H(x)(U)\subseteq U$ for all $x\in S$. $H$ is called \emph{irreducible} if $V$ and $\{0\}$ are the only $H$-invariant subspaces of $V$.}
\end{defin}

Note that the property of irreducibility depends on the field $F$. Since we only consider frameworks in the real vector space $\mathbb{R}^{d}$, the representation space of any linear representation in this paper is assumed to be a real vector space.

\begin{defin}
\emph{A linear representation $H:S\to GL(V)$ is said to be \emph{unitary} with respect to a given inner product $\langle v,w\rangle$ if \begin{displaymath}\langle H(x)(v),H(x)(w)\rangle=\langle v,w\rangle \quad\textrm{ for all } v,w\in V \textrm{ and all } x\in S\textrm{.}\end{displaymath} }
\end{defin}

\begin{remark}\label{unitary} \emph{A unitary representation has the property that the orthogonal complement of an invariant subspace is again invariant \cite{serre}.}
\end{remark}

\begin{defin}
\emph{Let $H:S\to GL(V)$ be a linear representation of a group $S$ and let $U$ be an invariant subspace of $V$. If for all $x\in S$, we restrict the automorphism $H(x)$ of $V$ to the subspace $U$, then we obtain a new linear representation $H^{(U)}$ of $S$ with representation space $U$. $H^{(U)}$ is said to be a \emph{subrepresentation} of $H$.}
\end{defin}

\begin{defin}
\emph{Let $H_{1}:S\to GL(V_{1})$ and $H_{2}:S\to GL(V_{2})$ be two linear representations of a group $S$. Then $H_{1}\oplus H_{2}:S\to GL(V_{1}\oplus V_{2})$ is the representation of $S$ which sends $x\in S$ to $H_{1}\oplus H_{2}(x)$, where $H_{1}\oplus H_{2}(x)\big((v_{1},v_{2})\big)=\big(H_{1}(x)(v_{1}),H_{2}(x)(v_{2})\big)$ for all $v_{1}\in V_{1}$ and $v_{2}\in V_{2}$.}
\end{defin}

\begin{defin}\emph{Let $S$ be a group and $F$ be a field. A \emph{matrix representation} of $S$ is a homomorphism $H$ from $S$ to $GL(n,F)$, where $GL(n,F)$ denotes the group of all invertible $n\times n$ matrices with entries in $F$.\\\indent Two matrix representations $H_{1}:S\to GL(n,F)$ and $H_{2}:S\to GL(n,F)$ are said to be \emph{equivalent} if there exists an invertible matrix $M$ such that $MH_{1}(x)M^{-1}=H_{2}(x)$ for all $x\in S$, in which case we write $H_{1}\backsimeq H_{2}$.}
\end{defin}

Let $S$ be a group, $V$ be an $n$-dimensional vector space over the field $F$, and $H:S\to GL(V)$ be a linear representation of $S$. Given a basis $B$ of $V$, we may associate a matrix representation $H_{B}:S\to GL(n,F)$ to $H$ by defining $H_{B}(x)$ to be the matrix that represents the automorphism $H(x)$ with respect to the basis $B$ for all $x\in S$. $H_{B}$ is then said to \emph{correspond} to $H$ with respect to $B$. Note that two matrix representations $H_{1}$ and $H_{2}$ correspond to equivalent linear representations if and only if $H_{1}\backsimeq H_{2}$.

\section{Block-diagonalization of the rigidity matrix}

\subsection{The internal and external representation}
\label{subsec:intext}

Given a graph $G$, a symmetry group $S$, and a homomorphism $\Phi:S\to \textrm{Aut}(G)$, we define two particular matrix representations of $S$, the external and the internal representation, both of which depend on $G$ and $\Phi$. These two representations play the key role in a symmetry-based rigidity analysis of a framework $(G,p)\in\mathscr{R}_{(G,S,\Phi)}$.\\\indent Note that our definitions of these representations are mathematically explicit definitions of the external and internal representation introduced in \cite{FGsymmax} and \cite{KG2}. Giving explicit definitions of these representations allows us to provide mathematical proofs for all the observations made in \cite{cfgsw, FGsymmax, KG2} and to extend these results so that they can also be applied to non-injective symmetric realizations in any dimension.

\begin{defin}
\label{inandexrep}
\emph{Let $G$ be a graph with $V(G)=\{v_{1},v_{2},\ldots,v_{n}\}$ and $E(G)=\{e_{1},e_{2},\ldots,e_{m}\}$, $S$ be a symmetry group in dimension $d$, and $\Phi$ be a homomorphism from $S$ to $\textrm{Aut}(G)$. For $x\in S$, let $M_{x}$ denote the orthogonal $d\times d$ matrix which represents $x$ with respect to the canonical basis of $\mathbb{R}^{d}$.\\\indent The \emph{external representation} of $S$ (with respect to $G$ and $\Phi$) is the matrix representation $H_{e}:S\to GL(dn,\mathbb{R})$ that sends $x\in S$ to the matrix $H_{e}(x)$ which is obtained from the transpose of the $n\times n$ permutation matrix corresponding to $\Phi(x)$ (with respect to the enumeration $V(G)=\{v_{1},v_{2},\ldots,v_{n}\}$) by replacing each 1 with the matrix $M_{x}$ and each 0 with a $d\times d$ zero-matrix.\\\indent The \emph{internal representation} of $S$ (with respect to $G$ and $\Phi$) is the matrix representation $H_{i}:S\to GL(m,\mathbb{R})$ that sends $x\in S$ to the transpose of the permutation matrix corresponding to the permutation of $E(G)$ (with respect to the enumeration $E(G)=\{e_{1},e_{2},\ldots,e_{m}\}$) which is induced by $\Phi(x)$.}
\end{defin}

\begin{remark}
\label{homrepre}
\emph{It is easy to verify that both the external representation $H_{e}$ and the internal representation $H_{i}$ of $S$ (with respect to $G$ and $\Phi$) are in fact matrix representations of the group $S$, provided that $\Phi$ is a homomorphism. If, however, $\Phi$ is not a homomorphism, then $H_{e}$ and $H_{i}$ are also not homomorphisms, in which case neither $H_{e}$ nor $H_{i}$ is a matrix representation of the group $S$.}
\end{remark}

\begin{examp}
\label{triangexam}
\emph{To illustrate the previous definition, let $K_{3}$ be the complete graph with $V(K_{3})=\{v_{1},v_{2},v_{3}\}$ and $E(K_{3})=\{e_{1},e_{2},e_{3}\}$, where $e_{1}=\{v_{1},v_{2}\}$, $e_{2}=\{v_{1},v_{3}\}$ and $e_{3}=\{v_{2},v_{3}\}$. Further, let $\mathcal{C}_{s}=\{Id,s\}$ be the symmetry group in dimension 2 with \begin{displaymath} M_{Id}= \left( \begin{array} {rr} 1 & 0 \\ 0 & 1 \end{array} \right)\textrm{ and } M_{s}= \left( \begin{array} {rr} -1 & 0 \\ 0 & 1 \end{array} \right)\textrm{,}\end{displaymath} and let $\Phi:\mathcal{C}_{s}\to \textrm{Aut}(K_{3})$ be the homomorphism defined by $\Phi(s)=(v_{1}\,v_{2})(v_{3})$. Then we have}

\begin{figure}[htp]
\begin{center}
\begin{tikzpicture}[very thick,scale=1]
\tikzstyle{every node}=[circle, draw=black, fill=white, inner sep=0pt, minimum width=5pt];
        \path (0,0) node (p1) [label = below left: $p_1$] {} ;
        \path (2,0) node (p2) [label = below right: $p_2$] {} ;
        \path (1,2) node (p3) [label = left: $p_3$] {} ;
        \draw (p1)  -- node [draw=white, below left=3pt] {$e_{1}$} (p2);
        \draw (p2) -- node [draw=white, right=2pt] {$e_{3}$} (p3);
        \draw (p3) -- node [draw=white, left=2pt] {$e_{2}$} (p1);
        \draw [dashed, thin] (1,-1) -- (1,3);
        \end{tikzpicture}
\end{center}
\caption{\emph{A framework $(K_{3},p)\in \mathscr{R}_{(K_{3},\mathcal{C}_{s},\Phi)}$.}}
\label{triang}
\end{figure}
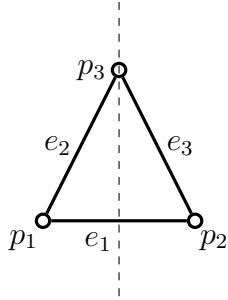

\begin{displaymath}
H_{e}(Id)= \left( \begin{array} {rr|rr|rr} 1 & 0 & 0 & 0 & 0 & 0\\ 0 & 1 & 0 & 0 & 0 & 0\\\hline 0 & 0 & 1 & 0 & 0 & 0\\0 & 0 & 0 & 1 & 0 & 0\\\hline 0 & 0 & 0 & 0 & 1 & 0\\0 & 0 & 0 & 0 & 0 & 1 \end{array} \right)\textrm{, } H_{e}(s)= \left( \begin{array} {rr|rr|rr} 0 & 0 & -1 & 0 & 0 & 0\\ 0 & 0 & 0 & 1 & 0 & 0\\\hline -1 & 0 & 0 & 0 & 0 & 0\\0 & 1 & 0 & 0 & 0 & 0\\\hline 0 & 0 & 0 & 0 & -1 & 0\\0 & 0 & 0 & 0 & 0 & 1\\\end{array} \right)\textrm{,}
\end{displaymath}
\begin{displaymath}
H_{i}(Id)= \left( \begin{array} {rrr} 1 & 0 & 0 \\ 0 & 1 & 0 \\ 0 & 0 & 1 \end{array} \right)\textrm{, } H_{i}(s)= \left( \begin{array} {rrr} 1 & 0 & 0\\ 0 & 0 & 1\\ 0 & 1 & 0\\\end{array} \right)\textrm{.}
\end{displaymath}
\end{examp}

For further examples, see \cite{KG2} or \cite{KG3}.

\subsection{The block-diagonalization}

In this section, we use the mathematically explicit definitions of the external and internal representation from the previous section to prove that the rigidity matrix of a symmetric framework can be transformed into a block-diagonalized form. Basic to this proof is Lemma \ref{replemma} which discloses the essential mathematical connection between the external and internal representation.

Recall from Section 2 that in the study of infinitesimal rigidity, we consider the equation \begin{displaymath}\mathbf{R}(G,p) u = z\textrm{,}\end{displaymath} where $\mathbf{R}(G,p)$ is the rigidity matrix of a framework $(G,p)$, $u\in \mathbb{R}^{d|V(G)|}$ is a column vector that represents an assignment of $d$-dimensional displacement vectors to the joints of $(G,p)$, and $z\in \mathbb{R}^{|E(G)|}$ is the column vector that represents the distortions in the bars of $(G,p)$ that are induced by $u$. The component of $z$ that corresponds to the edge $\{v_{i},v_{j}\}$ of $G$ is also known as the \emph{strain} induced on the bar $\{(v_{i},p_{i}),(v_{j},p_{j})\}$ by $u$.\\\indent Similarly, in the study of static rigidity, we consider the equation \begin{displaymath}\mathbf{R}(G,p)^{T} \omega = l\textrm{,}\end{displaymath} where the column vector $\omega\in \mathbb{R}^{|E(G)|}$ is a stress of $(G,p)$ and the column vector $l\in \mathbb{R}^{d|V(G)|}$ is the load on $(G,p)$ which is resolved by $\omega$.\\\indent Now, suppose $(G,p)$ is a symmetric framework in the set $\mathscr{R}_{(G,S,\Phi)}$, where $S$ is a symmetry group in dimension $d$ and $\Phi:S\to \textrm{Aut}(G)$ is a homomorphism. Then, using the notation of Definition \ref{inandexrep}, and assuming that the $i$th row of the rigidity matrix $\mathbf{R}(G,p)$ of $(G,p)$ corresponds to the edge $e_{i}$ of $G$, we have the following fundamental property of the external and internal representation of $S$ (with respect to $G$ and $\Phi$).

\begin{lemma}
\label{replemma}
Let $G$ be a graph, $S$ be a symmetry group, $\Phi$ be a homomorphism from $S$ to $\textrm{Aut}(G)$, and $p\in \bigcap_{x \in S} L_{x,\Phi}$.
\begin{itemize}
\item[(i)] If $\mathbf{R}(G,p) u = z$, then for all $x\in S$, we have $\mathbf{R}(G,p) H_{e}(x)u = H_{i}(x)z$;
\item[(ii)] if $\mathbf{R}(G,p)^{T} \omega = l$, then for all $x\in S$, we have $\mathbf{R}(G,p)^{T} H_{i}(x)\omega = H_{e}(x)l$.
\end{itemize}
\end{lemma}
\textbf{Proof.} $(i)$ Suppose $\mathbf{R}(G,p) u = z$. Fix $x\in S$ and let $M_{x}$ be the orthogonal matrix representing $x$ with respect to the canonical basis of $\mathbb{R}^{d}$. Also, let $\Phi(x)(v_{i})=v_{k}$ and $\Phi(x)(v_{j})=v_{l}$, and let $e_{f}=\{v_{i},v_{j}\}$ and $e_{h}=\{v_{k},v_{l}\}$. Then, since $p\in \bigcap_{x \in S} L_{x,\Phi}$, we have \begin{displaymath}M_{x}p_{i}=p_{k}\quad\textrm{ and }\quad M_{x}p_{j}=p_{l}\textrm{.}\end{displaymath} By the definition of $H_{i}(x)$, we have \begin{displaymath}\big(H_{i}(x)z\big)_{h}=(z)_{f}\textrm{.}\end{displaymath} Similarly, it follows from the definition of $H_{e}(x)$ that if $u\in \mathbb{R}^{dn}$ is replaced by $H_{e}(x)u$, then $u_{k}\in \mathbb{R}^{d}$ is replaced by $M_{x}u_{i}$ and $u_{l}\in \mathbb{R}^{d}$ by $M_{x}u_{j}$. By the definition of $\mathbf{R}(G,p)$, we have
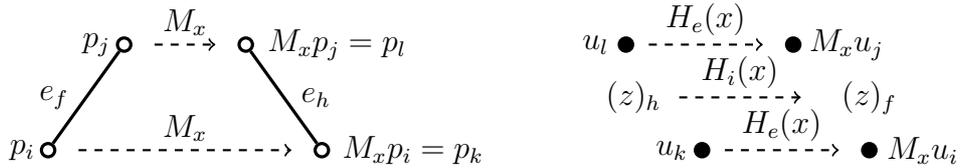
\begin{displaymath}
\big(\mathbf{R}(G,p) u\big)_{h}= (z)_{h}=(p_{k}-p_{l})\cdot u_{k}+(p_{l}-p_{k}) \cdot u_{l}\textrm{.}
\end{displaymath}

\begin{figure}[htp]
\begin{center}
\begin{tikzpicture}[very thick,scale=1]
\tikzstyle{every node}=[circle, draw=black, fill=white, inner sep=0pt, minimum width=5pt];
\path (0,0) node (pi) [label =  left: $p_i$] {} ;
\path (1,1.4) node (pj) [label =  left: $p_j$] {} ;
\node (pk) at (3.6,0) {};
\node [rectangle, draw=white ](label pk) at (4.8,0) {$M_xp_{i}= p_{k}$};
\node (pl) at (2.6,1.4) {};
\node [rectangle, draw=white](label pl) at (3.8,1.4) {$M_xp_{j}= p_{l}$};
\draw(pi)--node[rectangle, draw=white,left=6pt] {$e_f$}(pj);
\draw(pk)--node[rectangle, draw=white,right=6pt] {$e_h$}(pl);
\draw [->, thick, dashed] (0.4,0)--node[rectangle, draw=white,above=3pt] {$M_x$}(3.2,0);
\draw [->, thick, dashed] (1.4,1.4)--node[rectangle, draw=white,above=3pt] {$M_x$}(2.2,1.4);
\node[circle, fill=black] (uk) at (8.6,0) {};
\node[circle, fill=black] (ul) at (7.6,1.4) {};
\node [rectangle, draw=white ](label zh) at (7.7,0.7) {$ (z)_{h}$};
\node [rectangle, draw=white ](label uk) at (8.2,0) {$ u_{k}$};
\node [rectangle, draw=white ](label ul) at (7.2,1.4) {$u_{l}$};
\draw [->, thick, dashed] (8.3,0.7)--node[rectangle, draw=white,above=3pt] {$H_i(x)$}(9.9,0.7);
\draw [->, thick, dashed] (8.9,0)--node[rectangle, draw=white,above=3pt] {$H_e(x)$}(10.4,0);
\draw [->, thick, dashed] (7.9,1.4)--node[rectangle, draw=white,above =3pt] {$H_e(x)$}(9.4,1.4);
\node[circle, fill=black] (ukk) at (10.8,0) {};
\node[circle, fill=black] (ull) at (9.8,1.4) {};
\node [rectangle, draw=white ](label zf) at (10.8,0.7) {$ (z)_{f}$};
\node [rectangle, draw=white ](label uk) at (11.5,0) {$M_x u_{i}$};
\node [rectangle, draw=white ](label ul) at (10.5,1.4) {$M_x u_{j}$};
\end{tikzpicture}
\end{center}
\end{figure}

Therefore,
\begin{eqnarray}
\big(\mathbf{R}(G,p) H_{e}(x) u\big)_{h} & = & (p_{k}-p_{l}) \cdot M_{x}u_{i}+(p_{l}-p_{k}) \cdot M_{x}u_{j}\nonumber \\
& = & \big(M_{x}p_{i}-M_{x}p_{j}\big) \cdot M_{x}u_{i}+\big(M_{x}p_{j}-M_{x}p_{i}\big) \cdot M_{x}u_{j}\nonumber \\
& = & \big(M_{x}(p_{i}-p_{j})\big) \cdot M_{x}u_{i}+\big(M_{x}(p_{j}-p_{i})\big)\cdot M_{x}u_{j}\nonumber \\
& = & (p_{i}-p_{j})\cdot u_{i}+(p_{j}-p_{i})\cdot u_{j}\nonumber \\
& = & (z)_{f}\textrm{.}\nonumber
\end{eqnarray}
The penultimate equality sign is valid because the canonical inner product on $\mathbb{R}^{d}$ is invariant under the orthogonal transformation $x\in S$. This proves $(i)$.

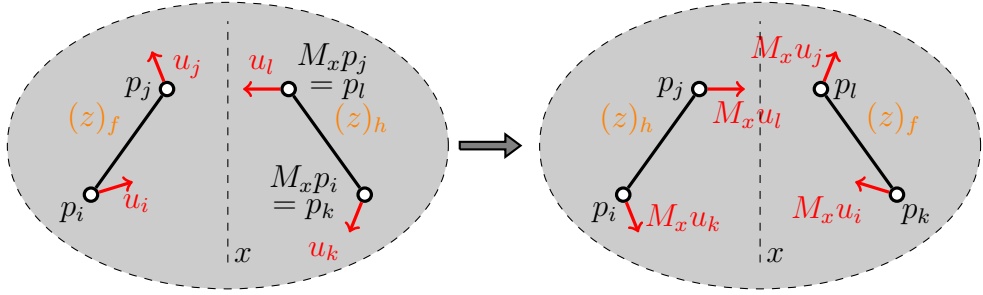
\begin{figure}[htp]
\begin{center}
\begin{tikzpicture}[very thick,scale=1]
\tikzstyle{every node}=[circle, draw=black, fill=white, inner sep=0pt, minimum width=5pt];
\filldraw[fill=black!20!white, draw=black, thin, dashed](1.8,0.65)ellipse(2.9cm and 1.9cm);
\path (0,0) node (pi) [label =  below left: $p_i$] {} ;
\path (1,1.4) node (pj) [label =  left: $p_j$] {} ;
\node (pk) at (3.6,0) {};
\node [rectangle, draw=black!20!white, fill=black!20!white](label pk) at (2.8,0.2) {$M_xp_{i}$};
\node [rectangle, draw=black!20!white, fill=black!20!white](label pk2) at (2.8,-0.2) {$= p_{k}$};
\node (pl) at (2.6,1.4) {};
\node [rectangle, draw=black!20!white, fill=black!20!white](label pl) at (3.18,1.8) {$M_xp_{j}$};
\node [rectangle, draw=black!20!white, fill=black!20!white](label pl) at (3.25,1.4) {$= p_{l}$};
\draw(pi)--node[color=orange,rectangle,draw=black!20!white, fill=black!20!white, above left=2pt] {$(z)_{f}$}(pj);
\draw(pk)--node[color=orange,rectangle,draw=black!20!white, fill=black!20!white,  above right=2pt] {$(z)_{h}$}(pl);
\draw [->, red] (pi) -- node[rectangle, draw=black!20!white, fill=black!20!white, below right=2pt] {$u_{i}$}(0.55,0.17);
\draw [->, red] (pj) -- node[rectangle, draw=black!20!white,fill=black!20!white, right=5pt] {$u_{j}$}(0.8,1.9);
\draw [->, red] (pk) -- (3.4,-0.5)node[rectangle, draw=black!20!white,fill=black!20!white, below left=3pt] {$u_{k}$};
\draw [->, red] (pl) -- node[rectangle,draw=black!20!white, fill=black!20!white, above=5pt] {$u_{l}$}(2,1.4);
\draw [dashed, thin] (1.8,-0.9)--(1.8,2.3);
\node [rectangle, draw=black!20!white, fill=black!20!white](label x) at (2,-0.8) {$x$};
\filldraw[fill=black!50!white, draw=black, thick]
    (4.85,0.6) -- (5.45,0.6) -- (5.45,0.5) -- (5.65,0.65) -- (5.45,0.8) -- (5.45,0.7) -- (4.85,0.7) -- (4.85,0.6);
\filldraw[fill=black!20!white, draw=black, thin, dashed](8.8,0.65)ellipse(2.9cm and 1.9cm);
\path (7,0) node (pii) [label =  below left: $p_i$] {} ;
\path (8,1.4) node (pjj) [label =  left: $p_j$] {} ;
\path (10.6,0) node (pkk) [label =  below right: $p_k$] {} ;
\path (9.6,1.4) node (pll) [label =  right: $p_l$] {} ;
\draw(pii)--node[color=orange,rectangle, draw=black!20!white,fill=black!20!white, above left=2pt] {$(z)_{h}$}(pjj);
\draw(pkk)--node[color=orange,rectangle, draw=black!20!white,fill=black!20!white,  above right=2pt] {$(z)_{f}$}(pll);
\draw [->, red] (pii) -- node[rectangle, draw=black!20!white, fill=black!20!white, right=4pt] {$M_xu_{k}$}(7.2,-0.5);
\draw [->, red] (pkk) -- node[rectangle, draw=black!20!white, fill=black!20!white, below left=3pt] {$M_xu_{i}$}(10.05,0.17);
\draw [->, red] (pjj) node[rectangle, draw=black!20!white, fill=black!20!white, below right=4pt] {$M_xu_{l}$}-- (8.6,1.4);
\draw [->, red] (pll) -- (9.8,1.9)node[rectangle, draw=black!20!white, fill=black!20!white, left=4pt] {$M_xu_{j}$};
\draw [dashed, thin] (8.8,-0.9)--(8.8,2.3);
\node [rectangle, draw=black!20!white, fill=black!20!white](label x) at (9,-0.8) {$x$};
\end{tikzpicture}
\end{center}
\caption{\emph{Illustration of the proof of Lemma \ref{replemma} (i) in the case where $x$ is a reflection.}}
\label{fig:Lemma_illustr}
\end{figure}

$(ii)$ Suppose $\mathbf{R}(G,p)^{T} \omega = l$. Fix $x\in S$ and let $\Phi(x)(v_{i})=v_{k}$. Then, since $p\in \bigcap_{x \in S} L_{x,\Phi}$, we have \begin{displaymath}M_{x}p_{i}=p_{k}\textrm{.}\end{displaymath} Let  $v_{i_{1}},v_{i_{2}},\ldots,v_{i_{j}}$ be the vertices in $V(G)$ that are adjacent to $v_{i}$, and let $e_{f_{t}}=\{v_{i},v_{i_{t}}\}$ for $t=1,2,\ldots,j$. Further, choose an enumeration of the $j$ vertices that are adjacent to $v_{k}$ in such a way that \begin{displaymath}M_{x}p_{i_{t}}= p_{k_{t}}\textrm{, }\end{displaymath} and let $e_{h_{t}}=\{v_{k},v_{k_{t}}\}$ for $t=1,2,\ldots,j$. For the vertex $v_{k}$, the equation $\mathbf{R}(G,p)^{T}\omega=l$ yields the vector-equation
\begin{equation}
\label{eqlemma}
(p_{k}-p_{k_{1}})(\omega)_{h_{1}}+\ldots +(p_{k}-p_{k_{j}})(\omega)_{h_{j}}=l_{k}\textrm{.}
\end{equation}

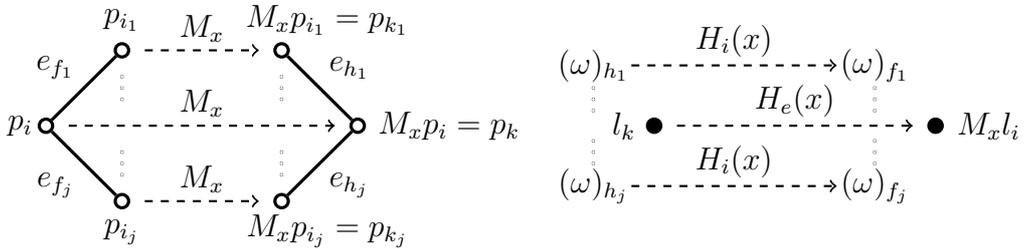
\begin{figure}[htp]
\begin{center}
\begin{tikzpicture}[very thick,scale=1]
\tikzstyle{every node}=[circle, draw=black, fill=white, inner sep=0pt, minimum width=5pt];
\path (0,0) node (pi) [label =  left: $p_i$] {} ;
\node (pk) at (4.1,0) {};
\node [rectangle, draw=white ](label pk) at (5.3,0) {$M_xp_{i}= p_{k}$};
\path (1,1) node (pi1) [label =  above: $p_{i_1}$] {} ;
\path (1,-1) node (pij) [label =  below: $p_{i_j}$] {} ;
\draw (pi)--node[rectangle, draw=white,above left=3pt] {$e_{f_1}$}(pi1);
\draw (pi)--node[rectangle, draw=white,below left=3pt] {$e_{f_j}$}(pij);
\node (pk1) at (3.1,1) {};
\node [rectangle, draw=white ](label pk1) at (3.7,1.4) {$M_xp_{i_1}= p_{k_1}$};
\node (pkj) at (3.1,-1) {};
\node [rectangle, draw=white ](label pkj) at (3.7,-1.4) {$M_xp_{i_j}= p_{k_j}$};
\draw (pk)--node[rectangle, draw=white,above right=3pt] {$e_{h_1}$}(pk1);
\draw (pk)--node[rectangle, draw=white,below right=3pt] {$e_{h_j}$}(pkj);
\draw [->, thick, dashed] (0.3,0)--node[rectangle, draw=white,above=3pt] {$M_x$}(3.8,0);
\draw [->, thick, dashed] (1.3,1)--node[rectangle, draw=white,above=3pt] {$M_x$}(2.8,1);
\draw [->, thick, dashed] (1.3,-1)--node[rectangle, draw=white,above=3pt] {$M_x$}(2.8,-1);
\filldraw[fill=black, draw=black]
    (1,0.65) circle (0.001cm);
\filldraw[fill=black, draw=black]
    (1,0.5) circle (0.001cm);
\filldraw[fill=black, draw=black]
    (1,0.35) circle (0.001cm);
\filldraw[fill=black, draw=black]
    (1,-0.65) circle (0.001cm);
\filldraw[fill=black, draw=black]
    (1,-0.5) circle (0.001cm);
\filldraw[fill=black, draw=black]
    (1,-0.35) circle (0.001cm);
\filldraw[fill=black, draw=black]
    (3.1,0.65) circle (0.001cm);
\filldraw[fill=black, draw=black]
    (3.1,0.5) circle (0.001cm);
\filldraw[fill=black, draw=black]
    (3.1,0.35) circle (0.001cm);
\filldraw[fill=black, draw=black]
    (3.1,-0.65) circle (0.001cm);
\filldraw[fill=black, draw=black]
    (3.1,-0.5) circle (0.001cm);
\filldraw[fill=black, draw=black]
    (3.1,-0.35) circle (0.001cm);
\node[circle, fill=black] (lk) at (8,0) {};
\node [rectangle, draw=white ](label lk) at (7.6,0) {$l_{k}$};
\node [rectangle, draw=white ](label oh1) at (7.2,0.8) {$(\omega)_{h_1}$};
\node [rectangle, draw=white ](label ohj) at (7.2,-0.8) {$(\omega)_{h_j}$};
\node[circle, fill=black] (lk) at (11.7,0) {};
\node [rectangle, draw=white ](label lk) at (12.4,0) {$M_xl_{i}$};
\node [rectangle, draw=white ](label oh1) at (10.9,0.8) {$(\omega)_{f_1}$};
\node [rectangle, draw=white ](label ohj) at (10.9,-0.8) {$(\omega)_{f_j}$};
\draw [->, thick, dashed] (8.3,0)--node[rectangle, draw=white,above=3pt] {$H_e(x)$}(11.4,0);
\draw [->, thick, dashed] (7.7,0.8)--node[rectangle, draw=white,above=3pt] {$H_i(x)$}(10.4,0.8);
\draw [->, thick, dashed] (7.7,-0.8)--node[rectangle, draw=white,above=3pt] {$H_i(x)$}(10.4,-0.8);
\filldraw[fill=black, draw=black]
    (7.2,0.5) circle (0.001cm);
\filldraw[fill=black, draw=black]
    (7.2,0.35) circle (0.001cm);
\filldraw[fill=black, draw=black]
    (7.2,0.2) circle (0.001cm);
\filldraw[fill=black, draw=black]
    (7.2,-0.5) circle (0.001cm);
\filldraw[fill=black, draw=black]
    (7.2,-0.35) circle (0.001cm);
\filldraw[fill=black, draw=black]
    (7.2,-0.2) circle (0.001cm);
\filldraw[fill=black, draw=black]
    (10.9,0.5) circle (0.001cm);
\filldraw[fill=black, draw=black]
    (10.9,0.35) circle (0.001cm);
\filldraw[fill=black, draw=black]
    (10.9,0.2) circle (0.001cm);
\filldraw[fill=black, draw=black]
    (10.9,-0.5) circle (0.001cm);
\filldraw[fill=black, draw=black]
    (10.9,-0.35) circle (0.001cm);
\filldraw[fill=black, draw=black]
    (10.9,-0.2) circle (0.001cm);
\end{tikzpicture}
\end{center}
\end{figure}

If $l\in \mathbb{R}^{dn}$ is replaced by $H_{e}(x)l$, then on the right-hand side of equation (\ref{eqlemma}), $l_{k}\in \mathbb{R}^{d}$ is replaced by $M_{x}l_{i}$ and if $\omega$ is replaced by $H_{i}(x)\omega$, then the left-hand side of equation (\ref{eqlemma}) is replaced by
\begin{eqnarray}
& &(p_{k}-p_{k_{1}})(\omega)_{f_{1}}+\ldots +(p_{k}-p_{k_{j}})(\omega)_{f_{j}}\nonumber \\
& = & \big(M_{x}p_{i}-M_{x}p_{i_{1}}\big)(\omega)_{f_{1}}+\ldots +\big(M_{x}p_{i}-M_{x}p_{i_{j}}\big)(\omega)_{f_{j}}\nonumber \\
& = & M_{x}\big((p_{i}-p_{i_{1}})(\omega)_{f_{1}}+\ldots +(p_{i}-p_{i_{j}})(\omega)_{f_{j}}\big)\nonumber \\
& = & M_{x}l_{i}\textrm{.}\nonumber
\end{eqnarray}
This completes the proof. $\square$

In the following, we again let $G$ be a graph with $V(G)=\{v_{1},v_{2},\ldots,v_{n}\}$ and $E(G)=\{e_{1},e_{2},\ldots,e_{m}\}$, $S$ be a symmetry group in dimension $d$, and $\Phi$ be a homomorphism from $S$ to $\textrm{Aut}(G)$.\\\indent Let $H_{e}$ be the external and $H_{i}$ be the internal representation of $S$ (with respect to $G$ and $\Phi$). Then we let $H'_{e}:S\to GL(\mathbb{R}^{dn})$ be the linear representation of $S$ that sends $x\in S$ to the automorphism $H'_{e}(x)$ which is represented by the matrix $H_{e}(x)$ with respect to the canonical basis of the $\mathbb{R}$-vector space $\mathbb{R}^{dn}$. Similarly, we let $H'_{i}:S\to GL(\mathbb{R}^{m})$ be the linear representation of $S$ that sends $x\in S$ to the automorphism $H'_{i}(x)$ which is represented by the matrix $H_{i}(x)$ with respect to the canonical basis of the $\mathbb{R}$-vector space $\mathbb{R}^{m}$. So, the external representation $H_{e}$ corresponds to the linear representation $H'_{e}$ with respect to the canonical basis of $\mathbb{R}^{dn}$ and the internal representation $H_{i}$ corresponds to the linear representation $H'_{i}$ with respect to the canonical basis of $\mathbb{R}^{m}$.\\\indent From group representation theory we know that every finite group has, up to equivalency, only finitely many irreducible linear representations and that every linear representation of such a group can be written uniquely, up to equivalency of the direct summands, as a direct sum of the irreducible linear representations of this group \cite{liebeck, serre}. So, let $S$ have $r$ pairwise non-equivalent irreducible linear representations $I_{1},I_{2},\ldots, I_{r}$ and let
\begin{equation}
\label{irrrep}
H'_{e}= \lambda_{1}I_{1}\oplus\ldots\oplus \lambda_{r}I_{r} \textrm{, where } \lambda_{1},\ldots,\lambda_{r}\in \mathbb{N}\cup {\{0\}} \textrm{.}
\end{equation}
For each $t=1,\ldots,r$, there exist $\lambda_{t}$ subspaces $\big(V_{e}^{(I_{t})}\big)_{1},\ldots, \big(V_{e}^{(I_{t})}\big)_{\lambda_{t}}$ of the $\mathbb{R}$-vector space $\mathbb{R}^{dn}$ which correspond to the $\lambda_{t}$ direct summands in (\ref{irrrep}), so that
\begin{equation}
\label{dirsumofvs}
\mathbb{R}^{dn}=V_{e}^{(I_{1})}\oplus \ldots \oplus V_{e}^{(I_{r})} \textrm{,}
\end{equation}
where
\begin{equation}
\label{dirsumofvs2}
V_{e}^{(I_{t})}= \big(V_{e}^{(I_{t})}\big)_{1}\oplus \ldots \oplus \big(V_{e}^{(I_{t})}\big)_{\lambda_{t}} \textrm{.}
\end{equation}
Let $\big(B_{e}^{(I_{t})}\big)_{1},\ldots, \big(B_{e}^{(I_{t})}\big)_{\lambda_{t}}$ be bases of the subspaces in (\ref{dirsumofvs2}). Then
\begin{equation}
\label{bases}
B_{e}^{(I_{t})}= \big(B_{e}^{(I_{t})}\big)_{1}\cup \ldots \cup \big(B_{e}^{(I_{t})}\big)_{\lambda_{t}}\nonumber
\end{equation}
is a basis of $V_{e}^{(I_{t})}$ and
\begin{equation}
\label{bases2}
B_{e}= B_{e}^{(I_{1})}\cup \ldots \cup B_{e}^{(I_{r})}
\end{equation}
is a basis of the $\mathbb{R}$-vector space $\mathbb{R}^{dn}$.

Consider now the matrix representation $\widetilde{H}_{e}$ that corresponds to the linear representation $H'_{e}$ with respect to the basis $B_{e}$. For $x\in S$, we have
\begin{equation}
\label{matr}
\widetilde{H}_{e}(x)=T_{e}^{-1}H_{e}(x)T_{e} \textrm{,}\nonumber
\end{equation}
where the $i$th column of $T_{e}$ is the coordinate vector of the $i$th basis vector of $B_{e}$ relative to the canonical basis, that is, $T_{e}$ is the matrix of the basis transformation from the canonical basis of the $\mathbb{R}$-vector space $\mathbb{R}^{dn}$ to the basis $B_{e}$. The column vectors of $\widetilde{H}_{e}(x)$ are the coordinates of the images of the basis vectors in $B_{e}$ under $H'_{e}(x)$ relative to the basis $B_{e}$. So, for each $x\in S$, the matrix $\widetilde{H}_{e}(x)$ has the same block form, namely
\begin{equation}
\label{Heblockform}\addtolength{\arraycolsep}{-0.8mm}
\widetilde{H}_{e}(x)=\left( \begin{array} {ccccccc} \big(A_{e}^{(I_{1})}\big)_{1}(x) &  &  &  &  & &\\ & \ddots &  &  &  & \mathbf{0} & \\  &  & \big(A_{e}^{(I_{1})}\big)_{\lambda_{1}}(x) &  &  &  & \\ &  &  & \ddots &  &  & \\  &  & &  & \big(A_{e}^{(I_{r})}\big)_{1}(x) & & \\ & \mathbf{0} &  &  &  &  \ddots & \\ &  & &  &  & & \big(A_{e}^{(I_{r})}\big)_{\lambda_{r}}(x)\end{array} \right)\textrm{.}\nonumber
\end{equation}
The block-matrix $\big(A_{e}^{(I_{t})}\big)_{j}(x)$ represents the restriction of the linear transformation $H'_{e}(x)$ to the subspace $\big(V_{e}^{(I_{t})}\big)_{j}$ with respect to the basis $\big(B_{e}^{(I_{t})}\big)_{j}$. Since for a given $t$, each of the subspaces $\big(V_{e}^{(I_{t})}\big)_{j}$, $j=1,\ldots,\lambda_{t}$, corresponds to the same irreducible linear representation $I_{t}$, we can choose the bases of the subspaces $\big(V_{e}^{(I_{t})}\big)_{j}$ in such a way that
\begin{equation}
\label{blocks}
\big(A_{e}^{(I_{t})}\big)_{1}(x)=\ldots = \big(A_{e}^{(I_{t})}\big)_{\lambda_{t}}(x)=:A_{e}^{(I_{t})}(x)\textrm{.}\nonumber
\end{equation}
In the following we assume that the basis $B_{e}$ is chosen in this way.

The above observations about the linear representation $H'_{e}$ of $S$ can be transferred analogously to the linear representation $H'_{i}$ of $S$. Let the direct sum decomposition of $H'_{i}$ be given by
\begin{equation}
\label{irrrepHi}
H'_{i}= \mu_{1}I_{1}\oplus\ldots\oplus \mu_{r}I_{r} \textrm{, where } \mu_{1},\ldots,\mu_{r}\in \mathbb{N}\cup {\{0\}} \textrm{.}
\end{equation}
For each $t=1,\ldots, r$, there exist $\mu_{t}$ subspaces $\big(V_{i}^{(I_{t})}\big)_{1},\ldots, \big(V_{i}^{(I_{t})}\big)_{\mu_{t}}$ of the $\mathbb{R}$-vector space $\mathbb{R}^{m}$ which correspond to the $\mu_{t}$ direct summands in (\ref{irrrepHi}), so that
\begin{equation}
\label{dirsumofvsHi}
\mathbb{R}^{m}=V_{i}^{(I_{1})}\oplus \ldots \oplus V_{i}^{(I_{r})} \textrm{,}
\end{equation}
where
\begin{equation}
\label{dirsumofvs2Hi}
V_{i}^{(I_{t})}= \big(V_{i}^{(I_{t})}\big)_{1}\oplus \ldots \oplus \big(V_{i}^{(I_{t})}\big)_{\mu_{t}} \textrm{.}
\end{equation}
Let $\big(B_{i}^{(I_{t})}\big)_{1},\ldots, \big(B_{i}^{(I_{t})}\big)_{\mu_{t}}$ be bases of the subspaces in (\ref{dirsumofvs2Hi}). Then
\begin{equation}
\label{basesHi}
B_{i}^{(I_{t})}= \big(B_{i}^{(I_{t})}\big)_{1}\cup \ldots \cup \big(B_{i}^{(I_{t})}\big)_{\mu_{t}}\nonumber
\end{equation}
is a basis of $V_{i}^{(I_{t})}$ and
\begin{equation}
\label{bases2Hi}
B_{i}= B_{i}^{(I_{1})}\cup \ldots \cup B_{i}^{(I_{r})}\nonumber
\end{equation}
is a basis of the $\mathbb{R}$-vector space $\mathbb{R}^{m}$.

Consider now the matrix representation $\widetilde{H}_{i}$ that corresponds to the linear representation $H'_{i}$ with respect to the basis $B_{i}$. Let $T_{i}$ be the matrix of the basis transformation from the canonical basis of the $\mathbb{R}$-vector space $\mathbb{R}^{m}$ to the basis $B_{i}$.  Then for $x\in S$, we have
\begin{equation}
\label{matrHi}
\widetilde{H}_{i}(x)=T_{i}^{-1}H_{i}(x)T_{i} \textrm{.}\nonumber
\end{equation}
So, the matrix $\widetilde{H}_{i}(x)$ has the same block form for each $x\in S$, namely
\begin{equation}
\label{Hiblockform}\addtolength{\arraycolsep}{-0.8mm}
\widetilde{H}_{i}(x)=\left( \begin{array} {ccccccc} \big(A_{i}^{(I_{1})}\big)_{1}(x) &  &  &  &  & &\\ & \ddots &  &  &  & \mathbf{0} & \\  &  & \big(A_{i}^{(I_{1})}\big)_{\mu_{1}}(x) &  &  &  & \\ &  &  & \ddots &  &  & \\  &  & &  & \big(A_{i}^{(I_{r})}\big)_{1}(x) & & \\ & \mathbf{0} &  &  &  &  \ddots & \\ &  & &  &  & & \big(A_{i}^{(I_{r})}\big)_{\mu_{r}}(x)\end{array} \right)\textrm{,}\nonumber
\end{equation}
and for each $t=1,2,\ldots,r$, we can choose the bases of the subspaces $\big(V_{i}^{(I_{t})}\big)_{j}$ in such a way that
\begin{equation}
\label{blocksHi}
\big(A_{i}^{(I_{t})}\big)_{1}(x)=\ldots = \big(A_{i}^{(I_{t})}\big)_{\mu_{t}}(x)=:A_{i}^{(I_{t})}(x)=A_{e}^{(I_{t})}(x)\textrm{.}\nonumber
\end{equation}
In the following we assume that $B_{i}$ is chosen in this way.

\begin{defin}
\label{usymmetric}
\emph{With the notation above, we say that a vector $v\in \mathbb{R}^{dn}$ is \emph{symmetric with respect to the irreducible linear representation $I_{t}$} of $S$ if $v\in V_{e}^{(I_{t})}$. Similarly, we say that a vector $w\in \mathbb{R}^{m}$ is \emph{symmetric with respect to the irreducible linear representation $I_{t}$} of $S$ if $w\in V_{i}^{(I_{t})}$.}
\end{defin}

We are now in the position to state the fundamental theorem for analyzing the rigidity properties of a symmetric framework using group representation theory.

\begin{theorem}
\label{fundthmrepth}
Let $G$ be a graph, $S$ be a symmetry group with pairwise non-equivalent irreducible linear representations $I_{1},\ldots,I_{r}$, $\Phi$ be a homomorphism from $S$ to $\textrm{Aut}(G)$, and $p\in \bigcap_{x \in S} L_{x,\Phi}$.
\begin{itemize}
\item[(i)] If $\mathbf{R}(G,p) u = z$ and $u$ is symmetric with respect to $I_{t}$, then $z$ is also symmetric with respect to $I_{t}$;
\item[(ii)] if $\mathbf{R}(G,p)^{T} \omega = l$ and $\omega$ is symmetric with respect to $I_{t}$, then $l$ is also symmetric with respect to $I_{t}$.
\end{itemize}
\end{theorem}
\textbf{Proof.} $(i)$ Suppose $S$ is a symmetry group in dimension $d$ and $G$ is a graph with $n$ vertices. Let $u\in \big(V_{e}^{(I_{t})}\big)_{j}$. By the direct sum decomposition of $V_{e}^{(I_{t})}$ in (\ref{dirsumofvs2}), the result follows if we can show that $z=\mathbf{R}(G,p) u\in V_{i}^{(I_{t})}$. By the decomposition of $\mathbb{R}^{|E(G)|}$ into direct summands in (\ref{dirsumofvs2Hi}), $z$ has a unique decomposition of the form
\begin{displaymath}
z=\sum_{\alpha=1}^{r}\sum_{\beta=1}^{\mu_{\alpha}}z_{\alpha,\beta} \textrm{, where } z_{\alpha,\beta}\in\big(V_{i}^{(I_{\alpha})}\big)_{\beta}\textrm{.}
\end{displaymath}
We now interpret $\mathbf{R}(G,p):\mathbb{R}^{dn}\to\mathbb{R}^{|E(G)|}$ as a linear transformation and for given $m$ and $k$, we define the projection map $\mathbf{R}_{m,k}$ corresponding to $\mathbf{R}(G,p)|_{\big(V_{e}^{(I_{t})}\big)_{j}}$ by
\begin{displaymath}
\mathbf{R}_{m,k}:\left\{\begin{array}{lll}\big(V_{e}^{(I_{t})}\big)_{j} & \to & \big(V_{i}^{(I_{m})}\big)_{k}\\
u & \mapsto & z_{m,k}\end{array}\right.\textrm{.}
\end{displaymath}
We need to show that for all $m\neq t$, $\mathbf{R}_{m,k}$ is the zero map. So, let $m\neq t$. Clearly, $\mathbf{R}_{m,k}$ is a linear transformation.\\\indent The image of $\mathbf{R}_{m,k}$ is an $H'_{i}$-invariant subspace of $\big(V_{i}^{(I_{m})}\big)_{k}$, as the following argument shows. Fix $x\in S$ and let $z'$ be in the image of $\mathbf{R}_{m,k}$, say $z'= \mathbf{R}_{m,k}(u')$. Then, by assumption, $H'_{e}(x)(u')\in \big(V_{e}^{(I_{t})}\big)_{j}$ and, by Lemma \ref{replemma} $(i)$, $H'_{i}(x)(z')$ is the image of $H'_{e}(x)(u')$ under $\mathbf{R}_{m,k}$. \\\indent
Since $I_{m}$ is an irreducible linear representation of $S$, $\big(V_{i}^{(I_{m})}\big)_{k}$ and $\{0\}$ are the only $H'_{i}$-invariant subspaces of $\big(V_{i}^{(I_{m})}\big)_{k}$. If the image of $\mathbf{R}_{m,k}$ is the null-space, then we are done, otherwise $\mathbf{R}_{m,k}$ is surjective.\\\indent Next, we show that the kernel of $\mathbf{R}_{m,k}$ is an $H'_{e}$-invariant subspace of $\big(V_{e}^{(I_{t})}\big)_{j}$. Fix $x\in S$ and let $u'$ be in the kernel of $\mathbf{R}_{m,k}$, that is, $\mathbf{R}_{m,k}(u')=0$. Then, again by Lemma \ref{replemma} $(i)$, the image of $H'_{e}(x)(u')$ under $\mathbf{R}_{m,k}$ is $H'_{i}(x)(0)=0$, and hence $H'_{e}(x)(u')$ is also in the kernel of $\mathbf{R}_{m,k}$.\\\indent Since $I_{t}$ is an irreducible linear representation of $S$, we either have $\textrm{ker }(\mathbf{R}_{m,k})=\big(V_{e}^{(I_{t})}\big)_{j}$, in which case we are done, or $\textrm{ker }(\mathbf{R}_{m,k})=\{0\}$, in which case $\mathbf{R}_{m,k}$ is injective.\\\indent So, assume $\mathbf{R}_{m,k}$ is bijective. Let the matrix that represents $\mathbf{R}_{m,k}$ with respect to the bases $\big(B_{e}^{(I_{t})}\big)_{j}$ and $\big(B_{i}^{(I_{m})}\big)_{k}$ be denoted by $\widetilde{\mathbf{R}}_{m,k}$. Then $\widetilde{\mathbf{R}}_{m,k}$ is an invertible matrix. Let $\tilde{u}$ be the coordinate vector of an element in $\big(V_{e}^{(I_{t})}\big)_{j}$ relative to the basis $\big(B_{e}^{(I_{t})}\big)_{j}$ and let $\tilde{z}$ be the coordinate vector of the image of $\tilde{u}$ under $\mathbf{R}_{m,k}$ relative to the basis $\big(B_{i}^{(I_{m})}\big)_{k}$. Then, by Lemma \ref{replemma} $(i)$, for any $x\in S$, we have
\begin{displaymath}
\widetilde{\mathbf{R}}_{m,k}\big(A_{e}^{(I_{t})}\big)_{j}(x)\tilde{u}=\big(A_{i}^{(I_{m})}\big)_{k}(x)\tilde{z}=
\big(A_{i}^{(I_{m})}\big)_{k}(x)\widetilde{\mathbf{R}}_{m,k}\tilde{u}\textrm{,}
\end{displaymath}
and hence also
\begin{displaymath}
\widetilde{\mathbf{R}}_{m,k}\big(A_{e}^{(I_{t})}\big)_{j}(x)=
\big(A_{i}^{(I_{m})}\big)_{k}(x)\widetilde{\mathbf{R}}_{m,k}\textrm{.}
\end{displaymath}
Therefore,
\begin{displaymath}
\widetilde{\mathbf{R}}_{m,k}\big(A_{e}^{(I_{t})}\big)_{j}(x)\widetilde{\mathbf{R}}^{-1}_{m,k}=
\big(A_{i}^{(I_{m})}\big)_{k}(x)=\big(A_{e}^{(I_{m})}\big)_{k}(x)\quad\textrm{ for all } x\in S\textrm{,}
\end{displaymath}
which says that $I_{t}$ and $I_{m}$ are equivalent representations, a contradiction. This completes the proof of part $(i)$.\\\indent With the help of Lemma \ref{replemma} $(ii)$, part $(ii)$ can be proved completely analogously to part $(i)$. $\square$

Theorem \ref{fundthmrepth} $(i)$ says that if $u\in \mathbb{R}^{dn}$ is an assignment of displacement vectors to the joints of a framework $(G,p)\in\mathscr{R}_{(G,S,\Phi)}$ and $u$ is symmetric with respect to $I_{t}$, then the strains induced on the bars of $(G,p)$ by $u$ must also be symmetric with respect to $I_{t}$. Similarly, Theorem \ref{fundthmrepth} $(ii)$ says that if $\omega$ is a resolution of an equilibrium load $l$ on $(G,p)\in\mathscr{R}_{(G,S,\Phi)}$ and $\omega$ is symmetric with respect to $I_{t}$, then $l$ must also be symmetric with respect to $I_{t}$.\\\indent An immediate consequence of Theorem \ref{fundthmrepth} is that the matrices $\mathbf{R}(G,p)$ and $\mathbf{R}(G,p)^{T}$ can be block-diagonalized in such a way that the original rigidity problems $\mathbf{R}(G,p) u = z$ and $\mathbf{R}(G,p)^{T} \omega = l$ are decomposed into subproblems, where each subproblem considers, respectively, the relationship between vectors $u$ and $z$ and vectors $\omega$ and $l$ that are symmetric with respect to the same irreducible linear representation $I_{t}$. This is specified in

\begin{cor}
\label{blockdiagon}
Let $G$ be a graph, $S$ be a symmetry group with pairwise non-equivalent irreducible linear representations $I_{1},\ldots,I_{r}$, $\Phi$ be a homomorphism from $S$ to $\textrm{Aut}(G)$, and $p\in \bigcap_{x \in S} L_{x,\Phi}$. Then the matrices $T_{i}^{-1}\mathbf{R}(G,p)T_{e}$ and $T_{e}^{-1}\mathbf{R}(G,p)^{T}T_{i}$ are block-diagonalized in such a way that there exists (at most) one submatrix block for each irreducible linear representation $I_{t}$ of $S$.
\end{cor}
\textbf{Proof.}
Suppose $\mathbf{R}(G,p) u = z$, and let $\tilde{u}$ be the coordinate vector of $u$ relative to the basis $B_{e}$ and $\tilde{z}$ be the coordinate vector of $z$ relative to the basis $B_{i}$. Further, let $\widetilde{\mathbf{R}}(G,p)$ be the matrix that represents the linear transformation $\mathbf{R}(G,p)$ with respect to the bases $B_{e}$ and $B_{i}$, that is, \begin{displaymath}\widetilde{\mathbf{R}}(G,p)=T_{i}^{-1}\mathbf{R}(G,p)T_{e}\textrm{.}\end{displaymath} Then, by changing coordinates relative to the canonical bases of $\mathbb{R}^{dn}$ and $\mathbb{R}^{m}$ into coordinates relative to the bases $B_{e}$ and $B_{i}$, the equation
\begin{displaymath}
\mathbf{R}(G,p) u = z
\end{displaymath}
is converted into the equation
\begin{displaymath}
\widetilde{\mathbf{R}}(G,p)\tilde{u}= \tilde{z}\textrm{.}
\end{displaymath}
By Theorem \ref{fundthmrepth} $(i)$, the matrix $\widetilde{\mathbf{R}}(G,p)$ is block-diagonalized in such a way that there exists (at most) one submatrix block for each irreducible linear representation $I_{t}$ of $S$ and the submatrix block corresponding to $I_{t}$ is a matrix of the size $\textrm{dim }\big(V_{i}^{(I_{t})}\big)\times \textrm{dim }\big(V_{e}^{(I_{t})}\big)$. In particular, a submatrix block can possibly be an `empty matrix' which has rows but no columns or alternatively columns but no rows.\\\indent
Similarly, if we denote $\tilde{\omega}$ to be the coordinate vector of $\omega$ relative to the basis $B_{i}$, $\tilde{l}$ to be the coordinate vector of $l$ relative to the basis $B_{e}$, and \begin{displaymath}\widetilde{\mathbf{R}}(G,p)^{T}=T_{e}^{-1}\mathbf{R}(G,p)^{T}T_{i}\textrm{,}\end{displaymath} then we may carry out the same changes of coordinates as above to convert the equation\begin{displaymath}\mathbf{R}(G,p)^{T} \omega = l\end{displaymath}  into the equation \begin{displaymath}\widetilde{\mathbf{R}}(G,p)^{T} \tilde{\omega}= \tilde{l}\textrm{.}\end{displaymath} By Theorem \ref{fundthmrepth} $(ii)$, the matrix $\widetilde{\mathbf{R}}(G,p)^{T}$ is again block-diagonalized in such a way that there exists (at most) one block for each $I_{t}$. $\square$

\begin{remark}
\emph{Note that the matrix $\widetilde{\mathbf{R}}(G,p)^{T}$ is equal to the transpose of the matrix $\widetilde{\mathbf{R}}(G,p)$ if and only if both of the matrices $T_{e}$ and $T_{i}$ are orthogonal matrices (i.e., $T_{e}^{-1}=T_{e}^{T}$ and $T_{i}^{-1}=T_{i}^{T}$) if and only if both $B_{e}$ and $B_{i}$ are orthonormal bases. Since the external and internal representation are both unitary representations (for all $x\in S$, $H_e(x)$ and $H_i(x)$ are orthogonal matrices), the invariant subspaces in (\ref{dirsumofvs}) and (\ref{dirsumofvsHi}) are mutually orthogonal (see \cite{faessler, serre}, for example). Thus, $B_{e}$ and $B_{i}$ can always be chosen to be orthonormal.}
\end{remark}

\begin{examp}
\label{blockmatrixexam}
\emph{Let $K_{3}$, $\mathcal{C}_{s}=\{Id,s\}$, and $\Phi$ be as in Example \ref{triangexam} and consider the framework $(K_{3},p)\in \mathscr{R}_{(K_{3},\mathcal{C}_{s},\Phi)}$ shown in Figures \ref{triang} and \ref{trisymloadres}, where \begin{displaymath}p_{1}=\left(\begin{array} {r} -1\\ 0 \end{array} \right)\textrm{, } p_{2}=\left(\begin{array} {r} 1\\ 0 \end{array} \right)\textrm{, and } p_{3}=\left(\begin{array} {r} 0\\ 2 \end{array} \right)\textrm{.}\end{displaymath} The rigidity matrix of $(K_{3},p)$ is given by
\setlength{\arraycolsep}{1.34pt}\begin{eqnarray*}\mathbf{R}(K_{3},p) &=& \left( \begin{array} {cccccc} (p_{1}-p_{2})_{1}& (p_{1}- p_{2})_{2} & (p_{2}- p_{1})_{1} & (p_{2}- p_{1})_{2} & 0 & 0\\ (p_{1}- p_{3})_{1}& (p_{1}- p_{3})_{2} & 0 & 0 & (p_{3}- p_{1})_{1} & (p_{3}- p_{1})_{2}\\ 0 & 0 & (p_{2}- p_{3})_{1}& (p_{2}- p_{3})_{2} & (p_{3}- p_{2})_{1} & (p_{3}- p_{2})_{2}\end{array} \right)\\&=&\setlength{\arraycolsep}{4pt}\left( \begin{array} {rrrrrr} -2 & 0 & 2 & 0 & 0 & 0\\ -1 & -2 & 0 & 0 & 1 & 2\\ 0 & 0 & 1 & -2 & -1 & 2 \end{array} \right)\textrm{.}\end{eqnarray*}
The symmetry group $\mathcal{C}_{s}$ has two non-equivalent irreducible linear representations both of which are of degree 1. In the Mulliken notation which is commonly used in chemistry and physics (see \cite{cotton}, for example), they are denoted by $A'$ and $A''$. $A'$ maps both $Id$ and $s$ to the identity transformation, whereas $A''$ maps $Id$ to the identity transformation and $s$ to the linear transformation $A''(s)$ which is defined by $A''(s)(x)=-x$ for all $x\in \mathbb{R}$. We have
\begin{displaymath}\mathbb{R}^{6}=V_{e}^{(A')}\oplus V_{e}^{(A'')}\end{displaymath} and \begin{displaymath}\mathbb{R}^{3}=V_{i}^{(A')}\oplus V_{i}^{(A'')}\textrm{.}\end{displaymath}
It is easy to see that the elements of the subspace $V_{e}^{(A')}$ of $\mathbb{R}^{6}$ are of the form
\begin{displaymath}\left( \begin{array}{r} u_{1}\\u_{2}\\-u_{1}\\u_{2}\\0\\u_{3}\end{array} \right)\textrm{, where } u_{1},u_{2},u_{3}\in\mathbb{R}\textrm{,}\end{displaymath}
(see Figure \ref{trisymloadres} (a)), so that an orthonormal basis $B_{e}^{(A')}$ of $V_{e}^{(A')}$ is given by
\begin{displaymath}B_{e}^{(A')}=\left\{\left( \begin{array}{c} \frac{1}{\sqrt{2}}\\0\\-\frac{1}{\sqrt{2}}\\0\\0\\0\end{array} \right),\left( \begin{array}{c} 0\\\frac{1}{\sqrt{2}}\\0\\\frac{1}{\sqrt{2}}\\0\\0\end{array} \right),\left( \begin{array}{r} 0\\0\\0\\0\\0\\1\end{array} \right)\right\}\textrm{.}\end{displaymath}}
\begin{figure}[ht]
\begin{center}
\begin{tikzpicture}[very thick,scale=1]
\tikzstyle{every node}=[circle, draw=black, fill=white, inner sep=0pt, minimum width=5pt];
        \path (0,0) node (p1) [label = below left: $p_1$] {} ;
        \path (2,0) node (p2) [label = below right: $p_2$] {} ;
        \path (1,2) node (p3) [label = left: $p_3$] {} ;
        \draw (p1)  -- node [draw=white, below left=3pt] {$e_{1}$} (p2);
        \draw (p2) -- node [draw=white, right=2pt] {$e_{3}$} (p3);
        \draw (p3) -- node [draw=white, left=2pt] {$e_{2}$} (p1);
        \draw [dashed, thin] (1,-1) -- (1,3);
        \draw [->, red] (p1) -- node [draw=white, left=8.5pt] {\setlength{\arraycolsep}{0.4pt}$\left(\begin{array}{r} u_{1}\\u_{2}\end{array}\right)$} (-0.7,0.7);
        \draw [->, red] (p2) -- node [draw=white, right=8.5pt] {\setlength{\arraycolsep}{0.4pt}$\left(\begin{array}{r} -u_{1}\\u_{2}\end{array}\right)$} (2.7,0.7);
        \draw [->, red] (p3) -- node [draw=white, right=3pt] {\setlength{\arraycolsep}{0.4pt}$\left(\begin{array}{c}0\\u_{3}\end{array}\right)$} (1,2.5);
        \node [draw=white, fill=white] (a) at (1,-1.4) {(a)};
        \path (0,-5) node (p1) [label = below left: $p_1$] {} ;
        \path (2,-5) node (p2) [label = below right: $p_2$] {} ;
        \path (1,-3) node (p3) [label = left: $p_3$] {} ;
        \draw (p1)  -- node [draw=white, below left=3pt] {$e_{1}$} (p2);
        \draw (p2) -- node [draw=white, right=2pt] {$e_{3}$} (p3);
        \draw (p3) -- node [draw=white, left=2pt] {$e_{2}$} (p1);
        \draw [dashed, thin] (1,-6) -- (1,-2);
        \draw [->, orange] (p1) -- node [draw=white, above left=3pt] {$z_{2}$}(0.25, -4.5);
        \draw [->, orange] (p3) -- (0.75, -3.5);
        \draw [->, orange] (p1) -- (-0.5,-5);
        \draw [->, orange] (p2) -- (2.5,-5);
        \draw [->, orange] (p2) -- node [draw=white, below left=12pt]{$z_{1}$} node [draw=white, above right=3pt] {$z_{2}$}(1.75,-4.5);
        \draw [->, orange] (p3) -- (1.25,-3.5);
        \node [draw=white, fill=white] (c) at (1,-6.4) {(c)};
        \end{tikzpicture}
        \hspace{0.4cm}
        \begin{tikzpicture}[very thick,scale=1]
\tikzstyle{every node}=[circle, draw=black, fill=white, inner sep=0pt, minimum width=5pt];
        \path (0,0) node (p1) [label = below left: $p_1$] {} ;
        \path (2,0) node (p2) [label = below right: $p_2$] {} ;
        \path (1,2) node (p3) [label = right: $p_3$] {} ;
        \draw (p1)  -- node [draw=white, below left=3pt] {$e_{1}$} (p2);
        \draw (p2) -- node [draw=white, right=2pt] {$e_{3}$} (p3);
        \draw (p3) -- node [draw=white, left=2pt] {$e_{2}$} (p1);
        \draw [dashed, thin] (1,-1) -- (1,3);
        \draw [->, red] (p1) -- node [draw=white, left=8.5pt] {\setlength{\arraycolsep}{0.4pt}$\left(\begin{array}{r} u_{1}\\u_{2}\end{array}\right)$} (-0.7,0.7);
        \draw [->, red] (p2)node [draw=white, right=11pt] {\setlength{\arraycolsep}{0.4pt}$\left(\begin{array}{r} u_{1}\\-u_{2}\end{array}\right)$} --  (1.3,-0.7);
        \draw [->, red] (p3) -- node [draw=white, left=7pt] {\setlength{\arraycolsep}{0.4pt}$\left(\begin{array}{c}u_{3}\\0\end{array}\right)$} (0.5,2);
        \node [draw=white, fill=white] (b) at (1,-1.4) {(b)};
        \path (0,-5) node (p1) [label = below left: $p_1$] {} ;
        \path (2,-5) node (p2) [label = below left: $p_2$] {} ;
        \path (1,-3) node (p3) [label = left: $p_3$] {} ;
        \draw (p1)  -- node [draw=white, below left=3pt] {$e_{1}$} (p2);
        \draw (p2) -- node [draw=white, right=2pt] {$e_{3}$} (p3);
        \draw (p3) -- node [draw=white, left=2pt] {$e_{2}$} (p1);
        \draw [dashed, thin] (1,-6) -- (1,-2);
        \draw [->, orange] (p1) -- node [draw=white, above left=3pt] {$z_{1}$}(0.25, -4.5);
        \draw [->, orange] (p3) -- (0.75, -3.5);
        \draw [->, orange] (p2) -- node [draw=white, left=19pt]{$0$} (2.25,-5.5)node [minimum width=2pt, draw=white, above = 18pt] {$-z_{1}$};
        \draw [->, orange] (p3) -- (0.75,-2.5);
        \node [draw=white, fill=white] (d) at (1,-6.4) {(d)};
        \end{tikzpicture}
\end{center}
\caption{\emph{\emph{(a, b)} Vectors of the $H_{e}'$-invariant subspaces $V_{e}^{(A')}$ \emph{(a)} and $V_{e}^{(A'')}$ \emph{(b)} of $\mathbb{R}^{6}$; \emph{(c, d)} vectors of the $H_{i}'$-invariant subspaces $V_{i}^{(A')}$ \emph{(c)} and $V_{i}^{(A'')}$ \emph{(d)} of $\mathbb{R}^{3}$.}}
\label{trisymloadres}
\end{figure}
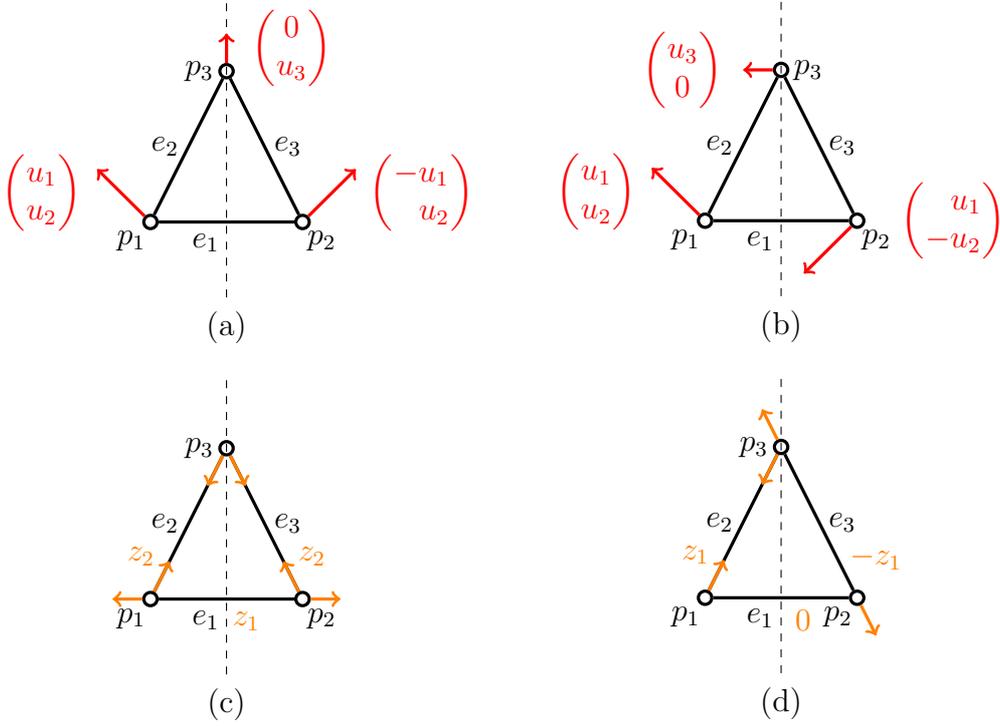
\emph{Similarly, the elements of the subspace $V_{e}^{(A'')}$ of $\mathbb{R}^{6}$ are of the form
\begin{displaymath}\left( \begin{array}{r} u_{1}\\u_{2}\\u_{1}\\-u_{2}\\u_{3}\\0\end{array} \right)\textrm{, where } u_{1},u_{2},u_{3}\in\mathbb{R}\textrm{,}\end{displaymath}
(see Figure \ref{trisymloadres} (b)), so that an orthonormal basis $B_{e}^{(A'')}$ of $V_{e}^{(A'')}$ is given by
\begin{displaymath}B_{e}^{(A'')}=\left\{\left( \begin{array}{c} \frac{1}{\sqrt{2}}\\0\\\frac{1}{\sqrt{2}}\\0\\0\\0\end{array} \right),\left( \begin{array}{c} 0\\\frac{1}{\sqrt{2}}\\0\\-\frac{1}{\sqrt{2}}\\0\\0\end{array} \right),\left( \begin{array}{r} 0\\0\\0\\0\\1\\0\end{array} \right)\right\}\textrm{.}\end{displaymath}
Orthonormal bases $B_{i}^{(A')}$ and $B_{i}^{(A'')}$ for the subspaces $V_{i}^{(A')}$ and $V_{i}^{(A'')}$ of $\mathbb{R}^{3}$ can be found analogously (see Figure \ref{trisymloadres} (c), (d)). We let \begin{displaymath}B_{i}^{(A')}=\left\{\left( \begin{array}{c} 1\\0\\0\end{array} \right),\left( \begin{array}{c} 0\\\frac{1}{\sqrt{2}}\\\frac{1}{\sqrt{2}}\end{array} \right)\right\}\textrm{.}\end{displaymath}
and
\begin{displaymath}B_{i}^{(A'')}=\left\{\left( \begin{array}{c} 0\\\frac{1}{\sqrt{2}}\\-\frac{1}{\sqrt{2}}\end{array} \right)\right\}\textrm{.}\end{displaymath}
Therefore, we have
\setlength{\arraycolsep}{2pt}\begin{displaymath}
T_{e}=\left( \begin{array} {cccccc} \frac{1}{\sqrt{2}} & 0 & 0 & \frac{1}{\sqrt{2}} & 0 & 0\\ 0 & \frac{1}{\sqrt{2}} & 0 & 0 & \frac{1}{\sqrt{2}} & 0\\ -\frac{1}{\sqrt{2}} & 0 & 0 & \frac{1}{\sqrt{2}} & 0 & 0\\ 0 & \frac{1}{\sqrt{2}} & 0 & 0 & -\frac{1}{\sqrt{2}} & 0 \\ 0 & 0 & 0 & 0 & 0 & 1\\ 0 & 0 & 1 & 0 & 0 & 0\end{array} \right)
\end{displaymath}
and
\setlength{\arraycolsep}{2pt}\begin{displaymath}
T_{i}=\left(\begin{array}{ccc} 1&0&0\\0&\frac{1}{\sqrt{2}}&\frac{1}{\sqrt{2}}\\0&\frac{1}{\sqrt{2}}&-\frac{1}{\sqrt{2}}\end{array} \right)\textrm{.}
\end{displaymath}
Thus,
\setlength{\arraycolsep}{2pt}\begin{displaymath}
\widetilde{\mathbf{R}}(K_{3},p)=T_{i}^{-1}\mathbf{R}(K_{3},p)T_{e}=\left( \begin{array} {ccc|ccc} -2\sqrt{2} & 0 & 0 & 0 & 0 & 0\\ -1 & -2 & 2\sqrt{2} & 0 & 0 & 0\\\hline 0 & 0 & 0 & -1 & -2 & \sqrt{2} \end{array} \right)
\end{displaymath}
and
\setlength{\arraycolsep}{2pt}\begin{displaymath}
\widetilde{\mathbf{R}}(K_{3},p)^{T}=T_{e}^{-1}\mathbf{R}(K_{3},p)^{T}T_{i}=\left( \begin{array} {cc|c} -2\sqrt{2} & -1 & 0\\ 0 & -2 & 0\\ 0 & 2\sqrt{2} & 0\\\hline 0 & 0 & -1\\ 0 & 0 & -2\\ 0 & 0 & \sqrt{2}\end{array} \right)\textrm{.}
\end{displaymath}
}
\end{examp}

\begin{remark}
\emph{In the previous example, we were able to find the invariant subspaces $V_{e}^{(A')}, V_{e}^{(A'')}$ of $\mathbb{R}^{6}$ and $V_{i}^{(A')}, V_{i}^{(A'')}$ of $\mathbb{R}^{3}$ by inspection because $\mathcal{C}_{s}$ is a small symmetry group with only two elements. This is of course generally not possible. There are, however, some standard methods and algorithms for finding the symmetry adapted bases $B_{e}$ and $B_{i}$ for any given symmetry group. Good sources for these methods are \cite{faessler, weeny}, for example.\\\indent As we will see in Section \ref{sec:maxwell}, knowledge of only the \emph{sizes} of the submatrix blocks that appear in the block-diagonalized rigidity matrices of a given symmetric framework allows us to gain significant insight into the rigidity properties of the framework. Since, with the aid of character theory, the sizes of these submatrix blocks can be determined very easily without explicitly finding the bases $B_{e}$ and $B_{i}$, there exist a number of applications of Corollary \ref{blockdiagon} (such as the symmetry-extended version of Maxwell's rule we will discuss in the following sections) that do not require finding the block-diagonalized rigidity matrices explicitly.
}
\end{remark}

\begin{remark}
\emph{The matrices $\mathbf{R}(G,p)^{T}\mathbf{R}(G,p)$ and $\mathbf{R}(G,p)\mathbf{R}(G,p)^{T}$ are also of interest in some areas of rigidity theory \cite{conwhit, KG2}. In structural engineering, these matrices are called the stiffness matrix and the flexibility matrix, respectively. It follows immediately from Corollary \ref{blockdiagon} that if $p\in \bigcap_{x \in S} L_{x,\Phi}$, then these matrices can also be block-diagonalized in such a way that there exists (at most) one block for each irreducible representation $I_{t}$ of $S$. In fact, it is easy to see that the matrices $T_{e}^{-1}\mathbf{R}(G,p)^{T}\mathbf{R}(G,p)T_{e}$ and $T_{i}^{-1}\mathbf{R}(G,p)\mathbf{R}(G,p)^{T}T_{i}$ have the desired block-form.}
\end{remark}

The fact that the rigidity matrix of a symmetric framework (as well as its transpose) can be block-diagonalized in the way described in Corollary \ref{blockdiagon} gives rise to many interesting results concerning the rigidity of symmetric frameworks \cite{cfgsw, FGsymmax, KG1, BS4, BS3, BS2}.\\\indent
Our goal for the remainder of this paper is to use Corollary \ref{blockdiagon} to establish a symmetry-extended version of Maxwell's rule that can be applied to (possibly non-injective) symmetric realizations in an arbitrary dimension $d$, and that contains the symmetry-extended version of Maxwell's rule stated in \cite{FGsymmax} as a special case.

\section{A symmetry-extended version of Maxwell's rule as an application}
\label{sec:maxwell}

\subsection{Maxwell's original rule}

If a framework is minimal infinitesimally (or statically) rigid, then it is said to be \emph{isostatic}. So, an isostatic framework is infinitesimally rigid and the removal of any bar results in a framework that is not infinitesimally rigid.\\\indent
Recall from Definition \ref{resolution} that a resolution of a load on a framework $(G,p)$ is also called a stress of $(G,p)$. A resolution of the zero-load is called a \emph{self-stress} of $(G,p)$. In other words, a self-stress is a linear dependence among the rows of the rigidity matrix of $(G,p)$.
If a framework does not have any non-zero self-stress, it is said to be \emph{independent}. Therefore, an isostatic framework is also characterized as  infinitesimally (or statically) rigid and independent \cite{gss, W1, W2}. In particular, the rows of the rigidity matrix of an isostatic framework $(G,p)$ form a basis for the space of equilibrium loads on $(G,p)$, provided that the points $p(v)$, $v\in V(G)$, span all of $\mathbb{R}^{d}$.

In 1864, Maxwell gave a necessary (but not sufficient) condition for a $2$- or $3$-dimensional framework $(G,p)$ to be isostatic \cite{bibmaxwell}. The $d$-dimensional version of this condition is given below. For additional necessary conditions, such as counts on all non-trivial subgraphs of $G$, see \cite{graver, gss, W1, W2}, for example.
In 1970, Laman provided sufficient conditions for `almost all' $2$-dimensional realizations of a given graph to be isostatic as well. However, there are well known problems in extending this result to higher dimensions \cite{graver, gss}.

\begin{theorem}[Maxwell's rule]
\label{Maxwell}
Let $(G,p)$ be a $d$-dimensional realization of a graph $G$ with $|V(G)|\geq d$. If $(G,p)$ is isostatic then \begin{displaymath}|E(G)|=d|V(G)|-\binom{d+1}{2}\textrm{.}\end{displaymath}
\end{theorem}

Let $(G,p)$ be a framework in $\mathbb{R}^d$ with the property that the points $p(v)$, $v\in V(G)$, span an affine subspace of $\mathbb{R}^d$ of dimension at least $d-1$, so that the space of infinitesimal rigid motions of $(G,p)$ has dimension $\binom{d+1}{2}$. Also, let the vector space of infinitesimal motions of $(G,p)$ be denoted by $I(p)$ and the vector space of self-stresses of $(G,p)$ be denoted by $\Omega(p)$. Then the equation in Maxwell's rule can be written in its extended form  as
\begin{displaymath}|E(G)|-d|V(G)|=\textrm{dim }\big(\Omega(p)\big)-\textrm{dim }\big(I(p)\big)\textrm{.}\end{displaymath} So, if $|E(G)|-\big(d|V(G)|-\binom{d+1}{2}\big)=k > 0$, then we can conclude that $(G,p)$ has at least $k$ linearly independent self-stresses and if $|E(G)|-\big(d|V(G)|-\binom{d+1}{2}\big)=-k < 0$, then $(G,p)$ has at least $k$ linearly independent infinitesimal flexes \cite{gss}.\\\indent
The advantage of Maxwell's rule is that it provides a purely combinatorial necessary condition for $(G,p)$ to be isostatic, and this condition can easily be verified since it only requires a simple count of the edges and vertices of $G$.

\subsection{The additional necessary conditions}

The symmetry-extended version of Maxwell's rule given in \cite{FGsymmax} provides further necessary conditions (in addition to Maxwell's original condition stated in Theorem \ref{Maxwell}) for a 2- or 3-dimensional symmetric framework with an injective configuration to be isostatic. Though the rule in \cite{FGsymmax} is a useful  tool for engineers and chemists to analyze the rigidity properties of symmetric structures in 2D and 3D, it is unsatisfactory from a mathematical point of view since it cannot be applied to frameworks in dimensions higher than 3, and since a complete mathematical proof of this result has not been provided. In the following sections, we aim to give a mathematical proof, based on the results of the previous sections, not only for the rule in \cite{FGsymmax}, but also for an extended rule that can be applied to a symmetric framework with a possibly non-injective configuration in an arbitrary dimension.\\\indent  In this section, we first develop all the necessary mathematical background that was omitted in \cite{FGsymmax}. This background consists of three major parts. First, we show that the subspaces $R$ and $T$ of all rotational and translational infinitesimal rigid motions of a given symmetric framework $(G,p)$ are invariant under the external representation $H'_e$ (Lemma \ref{Heinvar}), so that subrepresentations of $H'_e$ for the subspaces $R$ and $T$ can be defined. We then prove that the block-diagonalized form of the rigidity matrix of $(G,p)$ gives rise to additional necessary conditions for $(G,p)$ to be isostatic (Theorem \ref{Max}). The symmetry-extended version of Maxwell's rule is based on these conditions. Finally, we describe in detail how to determine the dimensions of the $H'_e$-invariant subspaces of $R$ and $T$. This is essential in applying the symmetry-extended version of Maxwell's rule to a given symmetric framework.\\\indent
Using some basic techniques from character theory, all of the results in this section combined will allow us to formulate the symmetry-extended version of Maxwell's rule given in \cite{FGsymmax} (as well as its extension to higher dimensions) as a mathematical theorem in Section \ref{subsec:rule}.\\\indent An alternate approach to proving the symmetry-extended version of Maxwell's rule given in \cite{FGsymmax} can be found in \cite{owen}.

In the following, we let $(G,p)$ be a symmetric framework in $\mathscr{R}_{(G,S,\Phi)}$, where $S$ is a non-trivial symmetry group in dimension $d$ and $\Phi:S\to \textrm{Aut}(G)$ is a homomorphism.\\\indent
In this section, we make the additional assumption that the points $p(v)$, $v\in V(G)$, span all of $\mathbb{R}^{d}$.\\\indent
Recall from Section 3 that we have the decomposition
\begin{equation}\label{Vedecomp1}
\mathbb{R}^{dn}=V_{e}^{(I_{1})}\oplus \ldots \oplus V_{e}^{(I_{r})}
\end{equation}
with
\begin{equation}
\label{Vedecomp}
V_{e}^{(I_{t})}= \big(V_{e}^{(I_{t})}\big)_{1}\oplus \ldots \oplus \big(V_{e}^{(I_{t})}\big)_{\lambda_{t}}
\end{equation}
of $\mathbb{R}^{dn}$ into $H'_{e}$-invariant subspaces.\\\indent
While the scalars $\lambda_{t}$ (as well as the subspaces that appear as direct summands in (\ref{Vedecomp1})) are uniquely determined in this decomposition, the subspaces that appear as direct summands in (\ref{Vedecomp}) are not \cite{serre}. In order to derive the desired symmetry-extended version of Maxwell's rule, the subspaces in (\ref{Vedecomp}) shall now be chosen appropriately.

Since the points $p(v)$, $v\in V(G)$, span all of $\mathbb{R}^{d}$, the subspace $N= \textrm{ker }\big(\mathbf{R}(K_{n},p)\big)$ of $\mathbb{R}^{dn}$, where $K_{n}$ is the complete graph on $V(G)$, is the space consisting of all infinitesimal rigid motions of $(G,p)$. This space can be written as the direct sum
\begin{displaymath}
N=T\oplus R\textrm{,}
\end{displaymath}
where $T$ is the space of all translational and $R$ is the space of all rotational infinitesimal rigid motions of $(G,p)$. More precisely, a basis of $T$ is given by $\{T_{j}|\, j=1,\ldots,d\}$, where for $j=1,\ldots,d$, $T_{j}:V(G)\to \mathbb{R}^{d}$ is the map that sends each $v\in V(G)$ to the $j$th canonical basis vector $e_{j}$ of $\mathbb{R}^{d}$, and a basis of $R$ is given by $\{R_{ij}|\,1\leq i<j\leq d\}$, where for $1\leq i<j\leq d$, $R_{ij}:V(G)\to \mathbb{R}^{d}$ is the map defined by $R_{ij}(v_{k})=(p_{k})_{i}e_{j}-(p_{k})_{j}e_{i}$ for all $k=1,\ldots,n$ \cite{W1}. Each of the maps $T_{j}$ and $R_{ij}$ is of course identified with a vector in $\mathbb{R}^{dn}$ (by using the order on $V(G)$).\\\indent
Note that in the context of static rigidity, $T$ is the space of all translational loads and $R$ is the space of all rotational loads on $(G,p)$.\\\indent
Using the notation of the previous paragraph we have the following result.

\begin{lemma}
\label{Heinvar}
For every dimension $d$, the subspaces $T$, $R$, and $N$ of $\mathbb{R}^{dn}$ are $H'_{e}$-invariant.
\end{lemma}
\textbf{Proof.} Fix a dimension $d$. We show first that $N=\textrm{ker }\big(\mathbf{R}(K_{n},p)\big)$ is $H'_{e}$-invariant. Since
$p\in \bigcap_{x \in S} L_{x,\Phi}$, it follows from Lemma \ref{replemma} that if
$\mathbf{R}(K_{n},p) u = z$, then for all $x\in S$, we have \begin{equation}\label{trneq} \mathbf{R}(K_{n},p) H_{e}(x)u = \widehat{H}_{i}(x)z\textrm{,}\end{equation} where $\widehat{H}_{i}$ is the internal representation of $S$ with respect to $K_{n}$ and $\Phi$. Let $u\in N$, i.e., $\mathbf{R}(K_{n},p) u = 0$. Then for any $x\in S$, we have
\begin{displaymath}\widehat{H}_{i}(x)\mathbf{R}(K_{n},p) u = \widehat{H}_{i}(x) 0 = 0\textrm{.}\end{displaymath} By (\ref{trneq}), we have  $\widehat{H}_{i}(x)\mathbf{R}(K_{n},p) u = \mathbf{R}(K_{n},p)H_{e}(x)u $, and hence
\begin{displaymath}\mathbf{R}(K_{n},p)H_{e}(x) u = 0\textrm{.}\end{displaymath} Thus, for all $x\in S$, $H_{e}(x) u \in \textrm{ker }\big(\mathbf{R}(K_{n},p)\big)$, which says that $N$ is $H'_{e}$-invariant.\\\indent Next, we show that $T$ is also $H'_{e}$-invariant. Let $x\in S$ and let, as usual, $M_{x}$ denote the orthogonal matrix that represents $x$ with respect to the canonical basis of $\mathbb{R}^{d}$. Then for $j=1,\ldots, d$, we have
\begin{displaymath}
H_{e}(x)T_{j}=\left(\begin{array}{c}M_{x}e_{j}\\\vdots\\M_{x}e_{j}\end{array}\right)=(M_{x})_{1j}T_{1}+\ldots + (M_{x})_{dj}T_{d}\textrm{.}
\end{displaymath}
It follows that $T$ is $H'_{e}$-invariant.\\\indent It remains to show that $R$ is $H'_{e}$-invariant. Since for all $x\in S$, $H_{e}(x)$ is an orthogonal matrix, $H'_{e}$ is a unitary representation (with respect to the canonical inner product on $\mathbb{R}^{dn}$). Therefore, the subrepresentation $H'^{(N)}_{e}$ of $H'_{e}$ with representation space $N$ is also unitary (with respect to the inner product obtained by restricting the canonical inner product on $\mathbb{R}^{dn}$ to $N$). So, by Remark \ref{unitary}, it suffices to show that $R$ is the orthogonal complement of $T$ in $N$.\\\indent Let $t$ be any element of $T$ and $r$ be any element of $R$. Then \begin{displaymath} t=\left(\begin{array}{c}w\\\vdots\\w\end{array}\right)\textrm{ for some } w\in \mathbb{R}^{d} \end{displaymath} and  \begin{displaymath}r=\left(\begin{array}{c}Vp_{1}\\\vdots\\Vp_{n}\end{array}\right)\textrm{ for some skew-symmetric matrix } V\textrm{.}\end{displaymath} Since the point $\sum_{i=1}^{n}p_{i}$ must be fixed by every symmetry operation $x\in S$, we may wlog define an origin so that $\sum_{i=1}^{n}p_{i}=0$. Then the inner product of $t$ and $r$ is given by
\begin{eqnarray} t \cdot r& =&\sum_{i=1}^{n}w^T Vp_{i}\nonumber\\
& = & w^T V\sum_{i=1}^{n}p_{i} = 0\textrm{.}\nonumber\end{eqnarray}
This gives the result. $\square$

Since, by Lemma \ref{Heinvar}, $N$ is an $H'_{e}$-invariant subspace of $\mathbb{R}^{dn}$, it follows from Maschke's Theorem (see \cite{liebeck, Mey, serre}, for example) that $N$ has an $H'_{e}$-invariant complement $Q$ in $\mathbb{R}^{dn}$.
We may therefore form the subrepresentation $H'^{(Q)}_{e}$ of $H'_{e}$ with representation space $Q$. Since $H'^{(Q)}_{e}$ is a direct sum of irreducible linear representations of $S$, say
\begin{equation}
\label{Qrep}
H'^{(Q)}_{e}=\kappa_{1}I_{1}\oplus\ldots\oplus \kappa_{r}I_{r}\textrm{, where }\kappa_{1},\ldots,\kappa_{r}\in \mathbb{N}\cup {\{0\}}\textrm{, }
\end{equation}
we obtain, analogously to (\ref{Vedecomp}), a decomposition of $Q$ of the form
\begin{displaymath}
Q=V_{Q}^{(I_{1})}\oplus \ldots \oplus V_{Q}^{(I_{r})} \textrm{,}
\end{displaymath}
where
\begin{equation}
\label{Qdecomp}
V_{Q}^{(I_{t})}= \big(V_{Q}^{(I_{t})}\big)_{1}\oplus \ldots \oplus \big(V_{Q}^{(I_{t})}\big)_{\kappa_{t}} \textrm{.}
\end{equation}
Similarly, since both $T$ and $R$ are also $H'_{e}$-invariant subspaces of $\mathbb{R}^{dn}$, we may form the subrepresentations $H'^{(T)}_{e}$ and $H'^{(R)}_{e}$ of $H'_{e}$ with respective representation spaces $T$ and $R$. This gives rise to a decomposition of $T$ of the form
\begin{displaymath}
T=V_{T}^{(I_{1})}\oplus \ldots \oplus V_{T}^{(I_{r})} \textrm{,}
\end{displaymath}
where
\begin{displaymath}
V_{T}^{(I_{t})}= \big(V_{T}^{(I_{t})}\big)_{1}\oplus \ldots \oplus \big(V_{T}^{(I_{t})}\big)_{\theta_{t}} \textrm{,}
\end{displaymath}
and to a decomposition of $R$ of the form
\begin{displaymath}
R=V_{R}^{(I_{1})}\oplus \ldots \oplus V_{R}^{(I_{r})} \textrm{,}
\end{displaymath}
where
\begin{displaymath}
V_{R}^{(I_{t})}= \big(V_{R}^{(I_{t})}\big)_{1}\oplus \ldots \oplus \big(V_{R}^{(I_{t})}\big)_{\rho_{t}} \textrm{.}
\end{displaymath}
We can now choose the decomposition in (\ref{Vedecomp}) in such a way that
\begin{equation}
\label{Veitdecomp}
V_{e}^{(I_{t})}= V_{Q}^{(I_{t})}\oplus V_{T}^{(I_{t})} \oplus V_{R}^{(I_{t})} \textrm{.}
\end{equation}
In the following we assume that the subspaces $\big(V_{e}^{(I_{t})}\big)_{j}$ are chosen in this way.

We are now in the position to derive the necessary conditions for $(G,p)\in\mathscr{R}_{(G,S,\Phi)}$ to be isostatic upon which the symmetry-extended version of Maxwell's rule is based.

\begin{theorem}
\label{Max}
Let $G$ be a graph, $S$ be a symmetry group in dimension $d$ with pairwise non-equivalent irreducible linear representations $I_{1},\ldots,I_{r}$, and $\Phi:S\to \textrm{Aut}(G)$ be a homomorphism. If $(G,p)$ is an isostatic framework in $\mathscr{R}_{(G,S,\Phi)}$ with the property that the points $p(v)$, $v\in V(G)$, span all of $\mathbb{R}^{d}$, then for $t=1,2,\ldots,r$, we have
\begin{equation}
\label{max}
\textrm{dim }\big(V_{Q}^{(I_{t})}\big)=\textrm{dim }\big(V_{i}^{(I_{t})}\big) \textrm{. }
\end{equation}
\end{theorem}
\textbf{Proof.} Suppose first that $\textrm{dim }\big(V_{Q}^{(I_{t})}\big)>\textrm{dim }\big(V_{i}^{(I_{t})}\big)$ for some $t$. In this case we give two separate arguments to show that $(G,p)$ is not isostatic, one that is based on infinitesimal rigidity and another one that is based on static rigidity. This will later allow us to obtain information about both kinematic and static rigidity properties of symmetric frameworks with the symmetry-extended version of Maxwell's rule.\\\indent It follows from Corollary \ref{blockdiagon} that there exists an element $u\neq 0$ in $V_{Q}^{(I_{t})}$ that lies in the kernel of the linear transformation which is represented by the matrix $\widetilde{\mathbf{R}}(G,p)$ with respect to the bases $B_{e}$ and $B_{i}$. In other words, $u$ is an infinitesimal flex of $(G,p)$ (which is symmetric with respect to $I_{t}$), and hence $(G,p)$ is not isostatic.\\\indent Alternatively, it follows from Corollary \ref{blockdiagon} that there exists an element $l$ in $V_{Q}^{(I_{t})}$ that does not lie in the image of the linear transformation which is represented by the matrix $\widetilde{\mathbf{R}}(G,p)^{T}$ with respect to the bases $B_{e}$ and $B_{i}$. This says that $l$ is an unresolvable equilibrium load on $(G,p)$ (which is symmetric with respect to $I_{t}$), so that we may again conclude that $(G,p)$ is not isostatic.\\\indent
Suppose now that $\textrm{dim }\big(V_{Q}^{(I_{t})}\big)<\textrm{dim }\big(V_{i}^{(I_{t})}\big)$ for some $t$. Then, analogously as above, there exists an element $\omega\neq 0$ in $V_{i}^{(I_{t})}$ that lies in the kernel of the linear transformation which is represented by the matrix $\widetilde{\mathbf{R}}(G,p)^{T}$ with respect to the bases $B_{e}$ and $B_{i}$. This says that $\omega$ is a non-zero self-stress of $(G,p)$ (which is symmetric with respect to $I_{t}$). So, it again follows that $(G,p)$ is not isostatic. $\square$

\begin{examp}\label{ex:rotandtranssp} \emph{Recall from Example \ref{blockmatrixexam} that for the framework $(K_3,p)\in \mathscr{R}_{(K_{3},\mathcal{C}_{s},\Phi)}$ shown in Figure \ref{fig:transandrot}, we have
\begin{eqnarray*}
\textrm{dim } \big(V_{e}^{(A')}\big) & = & 3\\
\textrm{dim }  \big(V_{i}^{(A')}\big)  &= & 2\\
\textrm{dim }  \big(V_{e}^{(A'')}\big)&  = & 3\\
\textrm{dim }  \big(V_{i}^{(A'')}\big) & = & 1 \textrm{.}
\end{eqnarray*}
It is easy to see that the 2-dimensional space $T$ of all translational infinitesimal rigid motions of $(K_3,p)$ can be written as the direct sum
\begin{displaymath}
T=V_{T}^{(A')} \oplus V_{T}^{(A'')}\textrm{,}
\end{displaymath}
where $V_{T}^{(A')}$ is the space of dimension 1 generated by the infinitesimal rigid motion shown in Figure \ref{fig:transandrot} (a), and $V_{T}^{(A'')}$ is the space of dimension 1 generated by the infinitesimal rigid motion shown in Figure \ref{fig:transandrot} (b). Moreover, the 1-dimensional space $R$ of rotational infinitesimal rigid motions of $(K_3,p)$ is clearly generated by the infinitesimal rigid motion shown in Figure \ref{fig:transandrot} (c), so that $R=V_{R}^{(A'')}$  and $\textrm{dim } \big(V_{R}^{(A')}\big)=0$.
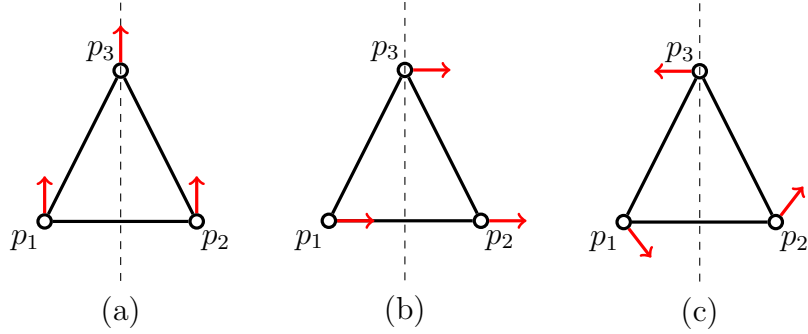
\begin{figure}[htp]
\begin{center}
\begin{tikzpicture}[very thick,scale=1]
\tikzstyle{every node}=[circle, draw=black, fill=white, inner sep=0pt, minimum width=5pt];
\path (0,0) node (p1) [label = below left: $p_1$] {} ;
        \path (2,0) node (p2) [label = below right: $p_2$] {} ;
        \path (1,2) node (p3) [label = above left: $p_3$] {} ;
        \draw (p1)  -- (p2);
        \draw (p2) --  (p3);
        \draw (p3) --  (p1);
        \draw [dashed, thin] (1,-1) -- (1,3);
\draw [->, red] (p1)--(0,0.6);
\draw [->, red] (p2)--(2,0.6);
\draw [->, red] (p3)--(1,2.6);
\node [draw=white, fill=white] (a) at (1,-1.2) {(a)};
\end{tikzpicture}
\hspace{0.5cm}
\begin{tikzpicture}[very thick,scale=1]
\tikzstyle{every node}=[circle, draw=black, fill=white, inner sep=0pt, minimum width=5pt];
\path (0,0) node (p1) [label = below left: $p_1$] {} ;
        \path (2,0) node (p2) [label = below right: $p_2$] {} ;
        \path (1,2) node (p3) [label = above left: $p_3$] {} ;
        \draw (p1)  -- (p2);
        \draw (p2) --  (p3);
        \draw (p3) --  (p1);
        \draw [dashed, thin] (1,-1) -- (1,3);
\draw [->, red] (p1)--(0.6,0);
\draw [->, red] (p2)--(2.6,0);
\draw [->, red] (p3)--(1.6,2);
\node [draw=white, fill=white] (b) at (1,-1.2) {(b)};
\end{tikzpicture}
\hspace{0.5cm}
\begin{tikzpicture}[very thick,scale=1]
\tikzstyle{every node}=[circle, draw=black, fill=white, inner sep=0pt, minimum width=5pt];
\path (0,0) node (p1) [label = below left: $p_1$] {} ;
        \path (2,0) node (p2) [label = below right: $p_2$] {} ;
        \path (1,2) node (p3) [label = above left: $p_3$] {} ;
        \draw (p1)  -- (p2);
        \draw (p2) --  (p3);
        \draw (p3) --  (p1);
        \draw [dashed, thin] (1,-1) -- (1,3);
\draw [->, red] (p1)--(0.36,-0.47);
\draw [->, red] (p2)--(2.36,0.47);
\draw [->, red] (p3)--(0.4,2);
\node [draw=white, fill=white] (c) at (1,-1.2) {(c)};
\end{tikzpicture}
\end{center}
\caption{\emph{\emph{(a)} A basis for the subspace $V_{T}^{(A')}$; \emph{(b)} a basis for the subspace $V_{T}^{(A'')}$; \emph{(c)} a basis for the subspace $R=V_{R}^{(A'')}$.}}
\label{fig:transandrot}
\end{figure}
It follows from equation (\ref{Veitdecomp}) that
\begin{equation}
\textrm{dim } \big(V_{Q}^{(A')}\big)=\textrm{dim } \big(V_{e}^{(A')}\big)-\textrm{dim } \big(V_{T}^{(A')}\big)-\textrm{dim } \big(V_{R}^{(A')}\big)=\textrm{dim } \big(V_{i}^{(A')}\big)=2\nonumber
\end{equation}
and
\begin{equation}
\textrm{dim } \big(V_{Q}^{(A'')}\big)=\textrm{dim } \big(V_{e}^{(A'')}\big)-\textrm{dim } \big(V_{T}^{(A'')}\big)-\textrm{dim } \big(V_{R}^{(A'')}\big)=\textrm{dim } \big(V_{i}^{(A'')}\big)=1\textrm{,}\nonumber
\end{equation}
so that the conditions (\ref{max}) in Theorem \ref{Max} are satisfied for the isostatic framework $(K_3,p)$.}
\end{examp}

In general, finding the dimensions of the subspaces $V_{Q}^{(I_{t})}$ and $V_{i}^{(I_{t})}$ by inspection is not as easy as in the previous example. In the following, we therefore describe a systematic method, based on techniques from character theory, for determining the dimensions of these subspaces, so that we can apply Theorem \ref{Max} to a symmetric framework with an arbitrary point group in any dimension. We begin by introducing the necessary vocabulary.

\begin{defin}
\emph{Let $A=(a_{ij})$ be an $n\times n$ square matrix. The \emph{trace} of $A$ is defined to be $Tr(A)=\sum_{i=1}^{n}a_{ii}$. }
\end{defin}

It is an important and well-known fact that the trace of a matrix is invariant under a similarity transformation \cite{cotton, Hall}. This gives rise to

\begin{defin}
\label{character}
\emph{Let $H:S\to GL(V)$ be a linear representation of a group $S$, $B$ be a basis of $V$, and $H_{B}$ be the matrix representation that corresponds to $H$ with respect to $B$. The \emph{character} $\chi(H)$ of $H$ is the function $\chi(H):S\to \mathbb{R}$ that sends $x\in S$ to $Tr\big(H_{B}(x)\big)$.\\\indent For a fixed enumeration $\{x_{1},\ldots,x_{k}\}$ of the elements of the group $S$, we will frequently also refer to the vector $\Big(Tr\big(H_{B}(x_{1})\big),\ldots,Tr\big(H_{B}(x_{k})\big)\Big)$ as the character of $H$.}
\end{defin}

In the following we need some well-known results from character theory which we summarize in

\begin{theorem}\cite{cotton, Hall, liebeck, serre}
\label{chartheorem}
Let $S$ be a group with $r$ pairwise non-equivalent irreducible linear representations $I_{1},\ldots, I_{r}$ and let $H:S\to GL(V)$ be a linear representation of $S$ with $H=\alpha_{1}I_{1}\oplus\ldots \oplus\alpha_{r}I_{r}$, where $\alpha_{t}\geq 0$ for all $t=1,\ldots, r$.
\begin{itemize}
\item[(i)] If $H=H_{1}\oplus H_{2}$ for some linear representations $H_{1}$ and $H_{2}$ of $S$, then $\chi(H)=\chi(H_{1})+\chi(H_{2})$;
\item[(ii)] $\chi(H)$ can be written uniquely as a linear combination of the characters $\chi(I_{1}),\ldots,\chi(I_{r})$ as \begin{displaymath}\chi(H)=\alpha_{1}\chi(I_{1})+\ldots +\alpha_{r}\chi(I_{r})\textrm{;}\end{displaymath}
\item[(iii)] For every $t=1,\ldots,r$, we have 
    \begin{displaymath}\alpha_{t} = \frac{1}{\|\chi(I_t)\|^2}\big( \chi(H) \cdot \chi(I_t)\big)\textrm{.}\end{displaymath}
\end{itemize}
\end{theorem}

We first explain how we can determine the dimensions of the subspaces $V_{i}^{(I_{t})}$ for all $t=1,\ldots,r$.\\\indent It follows from the direct sum decomposition of $H'_{i}$ in (\ref{irrrepHi}) that for $t=1,\ldots,r$, the dimension of $V_{i}^{(I_{t})}$ is the degree of $I_{t}$ multiplied by $\mu_{t}$. Since the degree of each irreducible linear representation $I_{t}$ can be read off from the character tables given in \cite{altherz, atkchil, Hall}, for example, we only need to determine the values of the $\mu_{t}$. This can easily be done by means of the formula given in Theorem \ref{chartheorem} $(iii)$, because the characters of the irreducible representations $I_{t}$ can simply be read off from the above-mentioned character tables and the character of $H'_{i}$ can be found by setting up the internal representation matrices $H_{i}(x)$, $x\in S$.\\\indent
Finding the dimensions of the subspaces $V_{Q}^{(I_{t})}$ for all $t=1,\ldots,r$ requires a little more work. It follows from (\ref{Veitdecomp}) that for $t=1,\ldots,r$, we have
\begin{displaymath}
\textrm{dim }\big(V_{Q}^{(I_{t})}\big)=\textrm{dim }\big(V_{e}^{(I_{t})}\big)- \textrm{dim }\big(V_{T}^{(I_{t})}\big)- \textrm{dim }\big(V_{R}^{(I_{t})}\big)\textrm{.}
\end{displaymath}
The dimensions of the subspaces $V_{e}^{(I_{t})}$ can be determined in the analogous way as the dimensions of the subspaces $V_{i}^{(I_{t})}$: for $t=1,\ldots,r$, the dimension of the subspace $V_{e}^{(I_{t})}$ is equal to the degree of $I_{t}$ multiplied by $\lambda_{t}$. Note that the values of the $\lambda_{t}$ in (\ref{Vedecomp}) can again easily be computed with the help of Theorem \ref{chartheorem} $(iii)$ since the character of $H'_{e}$ can be found by setting up the external representation matrices $H_{e}(x)$, $x\in S$. \\\indent
For $t=1,\ldots,r$, the dimension of the subspace $V_{T}^{(I_{t})}$ is the degree of $I_{t}$ multiplied by $\theta_{t}$ and the dimension of the subspace $V_{R}^{(I_{t})}$ is the degree of $I_{t}$ multiplied by $\rho_{t}$. So, in order to determine the dimensions of the subspaces $V_{T}^{(I_{t})}$ and $V_{R}^{(I_{t})}$ with the formula in Theorem \ref{chartheorem} $(iii)$, it only remains to determine the characters $\chi(H'^{(T)}_{e})$ and $\chi(H'^{(R)}_{e})$.\\\indent
We first show how to compute the character $\chi(H'^{(T)}_{e})$. It follows directly from the proof of Lemma \ref{Heinvar} that if $S$ is a symmetry group in dimension $d$ and $x\in S$, then the matrix that represents the linear transformation $H'^{(T)}_{e}(x)$ with respect to the basis $\{T_{1},\ldots,T_{d}\}$ is the orthogonal matrix $M_{x}$ that represents $x$ with respect to the canonical basis of $\mathbb{R}^{d}$. This says that for a fixed enumeration $\{x_{1},\ldots,x_{k}\}$ of the elements of $S$, we have
\begin{displaymath}
\chi(H'^{(T)}_{e})=\big(Tr(M_{x_{1}}),\ldots,Tr(M_{x_{k}})\big)\textrm{.}
\end{displaymath}
\indent For example, if $S$ is a symmetry group in dimension 2, then the component of $\chi(H'^{(T)}_{e})$ that corresponds to the identity in $S$ is equal to $2$, each component of $\chi(H'^{(T)}_{e})$ that corresponds to a rotation in $S$ about the origin by an angle of $\frac{2 \pi}{m}$ is equal to $2\cos\left(\frac{2 \pi}{m}\right)$, and each component of $\chi(H'^{(T)}_{e})$ that corresponds to a reflection in $S$ is equal to $0$.\\\indent
For a symmetry group in dimension 2 or 3, the explicit values of the components of $\chi(H'^{(T)}_{e})$ are summarized in \cite{cfgsw}.\\\indent
The character $\chi(H'^{(R)}_{e})$ can be determined similarly. As an example, we compute the character $\chi(H'^{(R)}_{e})$ in the case where $S$ is a symmetry group in dimension 2.
Every element of $S$ is then either the identity $Id$, a rotation $C_{m}$ about the origin by an angle of $\frac{2 \pi}{m}$, or a reflection $s$ in a line through the origin. Note that $R$ is a one-dimensional subspace of $\mathbb{R}^{2n}$ a basis of which is given by $\{R_{12}\}$. Let $C_{m}$ be a rotational symmetry operation of $(G,p)$ with
\begin{displaymath}
M_{C_{m}}=\left(\begin{array}{rr} \cos\left(\frac{2 \pi}{m}\right)&-\sin\left(\frac{2 \pi}{m}\right)\\\sin\left(\frac{2 \pi}{m}\right)&\cos\left(\frac{2 \pi}{m}\right)\end{array}\right)\textrm{.}
\end{displaymath}
Then, by using the definition of the external representation $H_{e}$ of $S$ (with respect to $G$ and $\Phi$) and the fact that $(G,p)\in\mathscr{R}_{(G,S,\Phi)}$, it is easy to verify that
\begin{displaymath}H_{e}(C_{m})R_{12}=R_{12}\textrm{.}
\end{displaymath}
Similarly, if $s$ is a reflectional symmetry operation of $(G,p)$ with
\begin{displaymath}M_{s}=\left(\begin{array}{rr} \cos\left(\theta\right)&\sin\left(\theta\right)\\\sin\left(\theta\right)&-\cos\left(\theta\right)\end{array}\right)\textrm{, }
\end{displaymath}
then it is again easy to verify that
\begin{displaymath}H_{e}(s)R_{12}=-R_{12}\textrm{.}
\end{displaymath}
It follows that the matrices which represent the linear transformations $H'^{(R)}_{e}(Id)$, $H'^{(R)}_{e}(C_{m})$, and $H'^{(R)}_{e}(s)$ with respect to the basis $\{R_{12}\}$ are the $1\times 1$ matrices (i.e., scalars) $1$, $1$, and $-1$, respectively. Therefore, if $d=2$, the character $\chi(H'^{(R)}_{e})$ is the vector defined as follows: each component of $\chi(H'^{(R)}_{e})$ that corresponds to the identity $Id\in S$ or a rotational symmetry operation $C_{m}\in S$ is equal to $1$, and each component of $\chi(H'^{(R)}_{e})$ that corresponds to a reflection $s\in S$ is equal to $-1$.\\\indent  Note that analogous calculations as above can easily be carried out for any symmetry group in dimension $d>2$ as well. For a symmetry group in dimension 2 or 3, the values of the components of $\chi(H'^{(R)}_{e})$ are again summarized in \cite{cfgsw}.

\begin{examp}\label{ex:dimcounts}
\emph{Let us apply the methods described above to the framework $(K_3,p)\in\mathscr{R}_{(K_{3},\mathcal{C}_{s},\Phi)}$ from Example \ref{ex:rotandtranssp}. From the representation matrices in Example \ref{triangexam} we immediately deduce that $\chi(H'_{e})=(6,0)$ and $\chi(H'_{i})=(3,1)$.  Therefore, if we let
\begin{eqnarray*}H'_{e}&=&\lambda_1 A' \oplus \lambda_2 A''\\
H'_{i}&=&\mu_1 A' \oplus \mu_2 A''\textrm{,}
\end{eqnarray*}
then, by the formula in Theorem \ref{chartheorem} $(iii)$, we have
\begin{eqnarray*}\lambda_1&=&\frac{1}{2}\big(6\cdot 1+0\cdot 1\big)=3\\
\lambda_2&=&\frac{1}{2}\big(6\cdot 1+0\cdot (-1)\big)=3\\
\mu_1&=&\frac{1}{2}\big(3\cdot 1+1\cdot 1\big)=2\\
\mu_2&=&\frac{1}{2}\big(3\cdot 1+1\cdot (-1)\big)=1\textrm{.}
\end{eqnarray*}
Further, for the characters $\chi(H'^{(T)}_{e})$ and $\chi(H'^{(R)}_{e})$, we have, as shown above, $\chi(H'^{(T)}_{e})=(2,0)$ and $\chi(H'^{(R)}_{e})=(1,-1)$. So, if we let
\begin{eqnarray*}H'^{(T)}_{e}&=&\theta_1 A' \oplus \theta_2 A''\\
H'^{(R)}_{e}&=&\rho_1 A' \oplus \rho_2 A''\textrm{,}
\end{eqnarray*}
then, again by the formula in Theorem \ref{chartheorem} $(iii)$, we obtain $\theta_1=1$, $\theta_2=1$, $\rho_1=0$, and $\rho_2=1$. Since both $A'$ and $A''$ are linear representations of degree 1, it follows that
\begin{eqnarray*}
\textrm{dim } \big(V_{Q}^{(A')}\big)&=&\textrm{dim } \big(V_{e}^{(A')}\big)-\textrm{dim } \big(V_{T}^{(A')}\big)-\textrm{dim } \big(V_{R}^{(A')}\big)\\&=&3-1-0\\&=&2\\
\textrm{dim } \big(V_{Q}^{(A'')}\big)&=&\textrm{dim } \big(V_{e}^{(A'')}\big)-\textrm{dim } \big(V_{T}^{(A'')}\big)-\textrm{dim } \big(V_{R}^{(A'')}\big)\\&=&3-1-1\\&=&1
\end{eqnarray*}
and
\begin{eqnarray*}
\textrm{dim } \big(V_{i}^{(A')}\big)&=&2\\
\textrm{dim } \big(V_{i}^{(A'')}\big)&=&1\textrm{.}
\end{eqnarray*}
}
\end{examp}

\subsection{The rule}
\label{subsec:rule}

Using  the mathematical background established in the previous section, we can now prove a symmetry-extended version of Maxwell's rule that can be applied to both injective and non-injective symmetric realizations in any dimension. Note that for dimensions 2 and 3, Theorem \ref{maxwellsrulewithchar} is a mathematically explicit formulation of the rule given in \cite{FGsymmax}.\\\indent The condition (\ref{maxchar}) in Theorem \ref{maxwellsrulewithchar} is obtained by combining all of the conditions in (\ref{max}) into a single equation using some basic techniques from character theory. This enables us to check the conditions in (\ref{max}) with very little computational effort, so that the symmetry-extended version of Maxwell's rule is in the same spirit as Maxwell's original rule in the sense that it only requires a simple count of joints and bars that are `fixed' by various symmetry operations.

From now on we will simplify our notation of the previous section by denoting the characters $\chi(H'_{e})$, $\chi(H'_{i})$, $\chi(H'^{(Q)}_{e})$, $\chi(H'^{(T)}_{e})$, and $\chi(H'^{(R)}_{e})$ by $X_{e}$, $X_{i}$, $X_{Q}$, $X_{T}$, and $X_{R}$, respectively.

\begin{theorem}[Symmetry-extended version of Maxwell's rule]
\label{maxwellsrulewithchar}
Let $G$ be a graph, $S$ be a symmetry group in dimension $d$ with pairwise non-equivalent irreducible linear representations $I_{1},\ldots,I_{r}$, and $\Phi:S\to \textrm{Aut}(G)$ be a homomorphism. If $(G,p)$ is an isostatic framework in $\mathscr{R}_{(G,S,\Phi)}$ with the property that the points $p(v)$, $v\in V(G)$, span all of $\mathbb{R}^{d}$, then we have
\begin{equation}
\label{maxchar}
X_{Q}= X_{i}\textrm{.}
\end{equation}
\end{theorem}
\textbf{Proof.} Suppose $X_{Q}\neq X_{i}$. Then, by Theorem \ref{chartheorem} $(ii)$ and equations (\ref{irrrepHi}) and (\ref{Qrep}), we have
\begin{displaymath}
\kappa_{1}\chi(I_{1})+\ldots +\kappa_{r}\chi(I_{r})\neq \mu_{1}\chi(I_{1})+\ldots +\mu_{r}\chi(I_{r})\textrm{, }
\end{displaymath}
which implies that $\kappa_{t}\neq \mu_{t}$ for some $t\in \{1,\ldots, r\}$. Therefore, $\textrm{dim }\big(V_{Q}^{(I_{t})}\big)\neq \textrm{dim }\big(V_{i}^{(I_{t})}\big)$. The result now follows from Theorem \ref{Max}. $\square$

So, by checking the condition (\ref{maxchar}) in Theorem \ref{maxwellsrulewithchar}, we implicitly check all the conditions in (\ref{max}). Since we have
\begin{displaymath}
H'_{e}=H'^{(Q)}_{e}\oplus H'^{(T)}_{e}\oplus H'^{(R)}_{e}\textrm{,}
\end{displaymath}
it follows from Theorem \ref{chartheorem} $(i)$ that
\begin{displaymath}
X_{Q}=X_{e}-X_{T}-X_{R}\textrm{.}
\end{displaymath}

Note that we have already shown how to compute each of the above characters in the previous section. In fact, for dimensions 2 and 3, the characters $X_{T}$ and $X_{R}$ can be read off from the tables in \cite{cfgsw}.
So, in order to check condition (\ref{maxchar}) for $d=2$ or $d=3$, it only remains to compute the characters $X_{e}$ and $X_{i}$.\\\indent So far, our method of determining $X_{e}$ and $X_{i}$ has been to set up the external and internal representation matrices $H_{e}(x)$ and $H_{i}(x)$ for all $x\in S$ and then to determine the traces of these matrices. In the following, we generalize the method of P. Fowler and S. Guest presented in \cite{FGsymmax} to determine the characters $X_{e}$ and $X_{i}$ without explicitly finding the external and internal representation of $S$. This will simplify significantly the calculations required to apply the symmetry-extended version of Maxwell's rule to a given framework.

\begin{defin}
\emph{Let $G$ be a graph with $V(G)=\{v_{1},\ldots,v_{n}\}$, $S$ be a symmetry group, $\Phi$ be a map from $S$ to $\textrm{Aut}(G)$, $(G,p)$ be a framework in $\mathscr{R}_{(G,S,\Phi)}$, and $x\in S$. A joint $(v_{i},p_{i})$ of $(G,p)$ is said to be \emph{fixed by $x$ with respect to $\Phi$} if $\Phi(x)(v_{i})=v_{i}$.\\\indent Similarly, a bar $\{(v_{i},p_{i}),(v_{j},p_{j})\}$ of $(G,p)$ is said to be \emph{fixed by $x$ with respect to $\Phi$} if $\Phi(x)\big(\{v_{i},v_{j}\}\big)=\{v_{i},v_{j}\}$.\\\indent The number of joints of $(G,p)$ that are fixed by $x$ with respect to $\Phi$ is denoted by $j_{\Phi(x)}$ and the number of bars of $(G,p)$ that are fixed by $x$ with respect to $\Phi$ is denoted by $b_{\Phi(x)}$.}
\end{defin}

Recall from Definition \ref{inandexrep} that for $x\in S$, the external representation matrix $H_{e}(x)$ is obtained from the transpose of the permutation matrix corresponding to $\Phi(x)$ by replacing each 1 with the $d\times d$ orthogonal matrix $M_{x}$ and each 0 with a $d\times d$ zero-matrix. Note that the transpose of the permutation matrix corresponding to $\Phi(x)$ has a 1 in the diagonal if and only if the corresponding joint of $(G,p)$ is fixed by $x$ with respect to $\Phi$. Therefore, a joint of $(G,p)$ can make a contribution to the trace of $H_{e}(x)$ only if it is fixed by $x$ with respect to $\Phi$. So, for a fixed enumeration $\{x_{1},\ldots,x_{k}\}$ of the elements of $S$, we have
\begin{eqnarray}
X_{e}&=&\Big(Tr\big(H_{e}(x_{1})\big),\ldots,Tr\big(H_{e}(x_{k})\big)\Big)\nonumber\\
&=&\big(j_{\Phi(x_{1})}Tr(M_{x_{1}}),\ldots,j_{\Phi(x_{k})}Tr(M_{x_{k}})\big)\nonumber\\
&=& X_{J}\times X_{T} \textrm{, }\nonumber
\end{eqnarray}
where $X_{J}=(j_{\Phi(x_{1})},\ldots,j_{\Phi(x_{k})})$
and $\times$ denotes componentwise multiplication.

Similarly, for $x\in S$, the internal representation matrix $H_{i}(x)$ has a 1 in the diagonal if and only if the corresponding bar of $(G,p)$ is fixed by $x$ with respect to $\Phi$. Thus,
\begin{displaymath}
X_{i}=(b_{\Phi(x_{1})},\ldots,b_{\Phi(x_{k})})\textrm{.}
\end{displaymath}

So, condition (\ref{maxchar}) in the symmetry-extended version of Maxwell's rule can be written as
\begin{equation}
\label{maxwellsequation}
X_{J}\times X_{T}-X_{T}-X_{R}=X_{i}\textrm{,}
\end{equation}
and each of the characters in (\ref{maxwellsequation}) can be determined with very little computational effort.

\begin{examp}\emph{The symmetry-extended version of Maxwell's rule, applied to the framework $(K_3,p)\in\mathscr{R}_{(K_{3},\mathcal{C}_{s},\Phi)}$ from Example \ref{ex:dimcounts}, yields the counts
\begin{eqnarray}
X_{J}&=&(j_{\Phi(Id)},j_{\Phi(s)})=(3,1)\nonumber\\
X_{T}&=&(2,0)\nonumber\\
X_{R}&=&(1,-1)\nonumber\\
X_{Q}&=& X_{J}\times X_{T}-X_{T}-X_{R}=(3,1)\nonumber\\
X_{i}&=&(b_{\Phi(Id)},b_{\Phi(s)})=(3,1)\textrm{.}\nonumber
\end{eqnarray}
Thus, condition (\ref{maxchar}) in Theorem \ref{maxwellsrulewithchar} is satisfied for the isostatic framework $(K_3,p)$.
}
\end{examp}

\begin{remark}\emph{Suppose the symmetry-extended version of Maxwell's rule detects that a framework $(G,p)\in\mathscr{R}_{(G,S,\Phi)}$ is not isostatic. Then we may use Theorem \ref{chartheorem} $(iii)$ and the proof of Theorem \ref{Max} to obtain information on the symmetry properties of self-stresses of $(G,p)$, infinitesimal flexes of $(G,p)$, and unresolvable equilibrium loads on $(G,p)$ in the following way.\\\indent Suppose for the framework $(G,p)$, we have $X_{Q}\neq X_{i}$. Using the formula in Theorem \ref{chartheorem} $(iii)$, we may then determine the values of the $\kappa_{t}$ and $\mu_{t}$ in (\ref{Qrep}) and (\ref{irrrepHi}) for all $t=1,\ldots,r$. By the proof of Theorem \ref{maxwellsrulewithchar}, there must exist $t\in \{1,\ldots, r\}$ such that $\kappa_{t}\neq \mu_{t}$.\\\indent It follows from the proof of Theorem \ref{Max} that if $\kappa_{t}> \mu_{t}$, say $\kappa_{t}- \mu_{t}=k>0$, then there exist $k$ linearly independent infinitesimal flexes of $(G,p)$ which are symmetric with respect to $I_{t}$, as well as $k$ unresolvable equilibrium loads on $(G,p)$ which are symmetric with respect to $I_{t}$.\\\indent Similarly, if $\kappa_{t}< \mu_{t}$, say $\mu_{t}- \kappa_{t}=k>0$, then there exist $k$ linearly independent self-stresses of $(G,p)$ which are symmetric with respect to $I_{t}$.}
\end{remark}

\subsection{Example and further remarks}

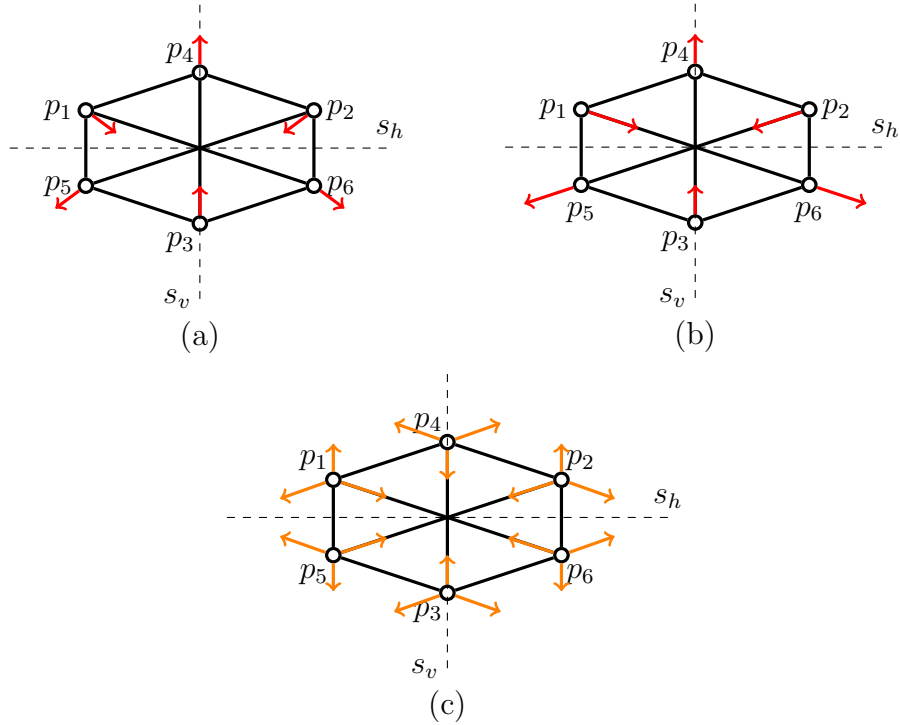
\begin{figure}[ht]
\begin{center}
\begin{tikzpicture}[very thick,scale=1]
\tikzstyle{every node}=[circle, draw=black, fill=white, inner sep=0pt, minimum width=5pt];
   \path (0,-0.5) node (p5) [label = left: $p_{5}$] {} ;
    \path (1.5,-1) node (p3) [label = below left: $p_{3}$] {} ;
    \path (3,-0.5) node (p6) [label = right: $p_{6}$] {} ;
   \path (0,0.5) node (p1) [label = left: $p_{1}$] {} ;
   \path (3,0.5) node (p2) [label = right: $p_{2}$] {} ;
   \path (1.5,1) node (p4) [label = above left: $p_{4}$] {} ;
   \draw (p1) -- (p4);
     \draw (p1) -- (p5);
     \draw (p1) -- (p6);
     \draw (p2) -- (p4);
     \draw (p2) -- (p5);
     \draw (p2) -- (p6);
     \draw (p3) -- (p4);
     \draw (p3) -- (p5);
     \draw (p3) -- (p6);
   \draw [dashed, thin] (1.5,-2) node [draw=white, left=2pt] {$s_{v}$}-- (1.5,2) ;
   \draw [dashed, thin] (-1,0) -- (4,0)node [draw=white, above=0.5pt] {$s_{h}$};
   \draw [->, red] (p5) -- (-0.4, -0.8);
   \draw [->, red] (p3) -- (1.5,-0.5);
   \draw [->, red] (p6) -- (3.4,-0.8);
   \draw [->, red] (p1) -- (0.4, 0.2);
   \draw [->, red] (p2) -- (2.6,0.2);
   \draw [->, red] (p4) -- (1.5,1.5);
      \node [draw=white, fill=white] (a) at (1.5,-2.5) {(a)};
\end{tikzpicture}
\hspace{1cm}
\begin{tikzpicture}[very thick,scale=1]
\tikzstyle{every node}=[circle, draw=black, fill=white, inner sep=0pt, minimum width=5pt];
\path (0,-0.5) node (p5) [label = below: $p_{5}$] {} ;
    \path (1.5,-1) node (p3) [label = below left: $p_{3}$] {} ;
    \path (3,-0.5) node (p6) [label = below: $p_{6}$] {} ;
   \path (0,0.5) node (p1) [label = left: $p_{1}$] {} ;
   \path (3,0.5) node (p2) [label = right: $p_{2}$] {} ;
   \path (1.5,1) node (p4) [label = above left: $p_{4}$] {} ;
   \draw (p1) -- (p4);
     \draw (p1) -- (p5);
     \draw (p1) -- (p6);
     \draw (p2) -- (p4);
     \draw (p2) -- (p5);
     \draw (p2) -- (p6);
     \draw (p3) -- (p4);
     \draw (p3) -- (p5);
     \draw (p3) -- (p6);
  \draw [dashed, thin] (1.5,-2) node [draw=white, left=2pt] {$s_{v}$}-- (1.5,2) ;
   \draw [dashed, thin] (-1,0) -- (4,0)node [draw=white, above=0.5pt] {$s_{h}$};
   \draw [->, red] (p5) -- (-0.75, -0.75);
   \draw [->, red] (p3) -- (1.5,-0.5);
   \draw [->, red] (p6) -- (3.75,-0.75);
   \draw [->, red] (p1) -- (0.75, 0.25);
   \draw [->, red] (p2) -- (2.25,0.25);
   \draw [->, red] (p4) -- (1.5,1.5);
      \node [draw=white, fill=white] (a) at (1.5,-2.5) {(b)};
\end{tikzpicture}
\hspace{1cm}
\begin{tikzpicture}[very thick,scale=1]
\tikzstyle{every node}=[circle, draw=black, fill=white, inner sep=0pt, minimum width=5pt];
\path (0,-0.5) node (p5) [label = below left: $p_{5}$] {} ;
    \path (1.5,-1) node (p3) [label = below left: $p_{3}$] {} ;
    \path (3,-0.5) node (p6) [label = below right: $p_{6}$] {} ;
   \path (0,0.5) node (p1) [label = above left: $p_{1}$] {} ;
   \path (3,0.5) node (p2) [label = above right: $p_{2}$] {} ;
   \path (1.5,1) node (p4) [label = above left: $p_{4}$] {} ;
   \draw (p1) -- (p4);
     \draw (p1) -- (p5);
     \draw (p1) -- (p6);
     \draw (p2) -- (p4);
     \draw (p2) -- (p5);
     \draw (p2) -- (p6);
     \draw (p3) -- (p4);
     \draw (p3) -- (p5);
     \draw (p3) -- (p6);
   \draw [dashed, thin] (1.5,-2) node [draw=white, left=2pt] {$s_{v}$}-- (1.5,2);
   \draw [dashed, thin] (-1.4,0) -- (4.4,0)node [draw=white, above=0.5pt] {$s_{h}$};
   \draw [->, orange] (p3) -- (1.5,-0.5);
   \draw [->, orange] (p4) -- (1.5,0.5);
   \draw [->, orange] (p4) -- (2.2, 1.25);
   \draw [->, orange] (p4) -- (0.8, 1.25);
   \draw [->, orange] (p3) -- (0.8, -1.25);
   \draw [->, orange] (p3) -- (2.2, -1.25);
   \draw [->, orange] (p1) -- (-0.7, 0.25);
   \draw [->, orange] (p5) -- (-0.7, -0.25);
   \draw [->, orange] (p2) -- (3.7, 0.25);
   \draw [->, orange] (p6) -- (3.7, -0.25);
   \draw [->, orange] (p1) -- (0.7, 0.271);
   \draw [->, orange] (p5) -- (0.7, -0.271);
   \draw [->, orange] (p6) -- (2.3, -0.271);
   \draw [->, orange] (p2) -- (2.3, 0.271);
   \draw [->, orange] (p1) -- (0, 0.978);
   \draw [->, orange] (p5) -- (0, -0.978);
   \draw [->, orange] (p2) -- (3, 0.978);
   \draw [->, orange] (p6) -- (3, -0.978);
      \node [draw=white, fill=white] (a) at (1.5,-2.5) {(c)};
\end{tikzpicture}
\end{center}
\caption{\emph{\emph{(a)} An infinitesimal flex of $(K_{3,3},p)\in\mathscr{R}_{(K_{3,3},\mathcal{C}_{2v},\Phi)}$ which is symmetric with respect to $B_{2}$ (the displacement vector at each joint of $(K_{3,3},p)$ remains unchanged under $Id$ and $s_{v}$ and is reversed under $C_{2}$ and $s_{h}$). \emph{(b)} An unresolvable equilibrium load on $(K_{3,3},p)$ which is symmetric with respect to $B_{2}$. \emph{(c)} A self-stress of $(K_{3,3},p)$ which is symmetric with respect to $A_{1}$ (the tensions and compressions in the bars of $(K_{3,3},p)$ remain unchanged under all symmetry operations in $\mathcal{C}_{2v}$).}}
\label{symfle}
\end{figure}

To illustrate how the symmetry-extended version of Maxwell's rule can give a significantly improved insight into the rigidity properties of a symmetric framework in comparison to Maxwell's original rule, we consider the framework $(K_{3,3},p)\in\mathscr{R}_{(K_{3,3},\mathcal{C}_{2v},\Phi)}$ depicted in Figure \ref{symfle}. $K_{3,3}$ is the complete bipartite graph from Example \ref{K33ex}, the symmetry group $\mathcal{C}_{2v}$ consists of the identity $Id$, the half-turn $C_{2}$, and the two reflections $s_{h}$ and $s_{v}$ whose mirror lines are the $x$-axis and the $y$-axis, respectively, and $\Phi: \mathcal{C}_{2v}\to \textrm{Aut}(K_{3,3})$ is defined by
\begin{eqnarray} \Phi(Id) & = & id\nonumber\\ \Phi(C_{2}) & = & (v_{1}\,v_{6})(v_{2}\, v_{5})(v_{3}\,v_{4})\nonumber\\ \Phi(s_{h}) & = & (v_{1}\,v_{5})(v_{2}\, v_{6})(v_{3}\,v_{4})\nonumber\\ \Phi(s_{v}) & = & (v_{1}\,v_{2})(v_{5}\,v_{6})(v_{3}) (v_{4})
\nonumber\textrm{.}\end{eqnarray}
The symmetry group $\mathcal{C}_{2v}$ has four non-equivalent irreducible linear representations each of which is of degree 1. In the Mulliken notation (\cite{cotton}), they are denoted by $A_{1}$, $A_{2}$, $B_{1}$, and $B_{2}$. The following table shows the characters of these representations:
\begin{table}[htp]
\begin{center}
\begin{tabular}{r||r|r|r|r}
$\mathcal{C}_{2v}$   &   $Id$  &  $C_{2}$   & $s_{h}$   &    $s_{v}$\\\hline\hline
$A_{1}$  &    1  &  1 &  1  &  1\\\hline
$A_{2}$  &    1  &  1 &  -1  &  -1\\\hline
$B_{1}$  &    1  &  -1 &  1  &  -1\\\hline
$B_{2}$  &    1  &  -1 &  -1  &  1\\
\end{tabular}
\end{center}
\end{table}

\noindent We have
\begin{eqnarray}
X_{J}&=&(j_{\Phi(Id)},j_{\Phi(C_{2})},j_{\Phi(s_{h})},j_{\Phi(s_{v})})=(6,0,0,2)\nonumber\\
X_{T}&=&(2,-2,0,0)\nonumber\\
X_{R}&=&(1,1,-1,-1)\nonumber\\
X_{Q}&=& X_{J}\times X_{T}-X_{T}-X_{R}=(9,1,1,1)\nonumber\\
X_{i}&=&(b_{\Phi(Id)},b_{\Phi(C_{2})},b_{\Phi(s_{h})},b_{\Phi(s_{v})})=(9,3,3,1)\textrm{.}\nonumber
\end{eqnarray}
Since $X_{Q}\neq X_{i}$, we may conclude that $(K_{3,3},p)$ is not isostatic. Note that Maxwell's original rule would not have detected this because $|E(K_{3,3})|=2|V(K_{3,3})|-3=9$.\\\indent With the help of the formula in Theorem \ref{chartheorem} $(iii)$ we obtain
\begin{eqnarray}
X_{Q}&=& 3A_{1}+2A_{2}+2B_{1}+ 2B_{2}\qquad\textrm{ and} \nonumber\\
X_{i}&=& 4A_{1}+2A_{2}+2B_{1}+ B_{2}\textrm{, }\nonumber
\end{eqnarray}
which implies that $(K_{3,3},p)$ has a non-zero self-stress which is symmetric with respect to $A_{1}$ and an infinitesimal flex (as well as an unresolvable equilibrium load) which is symmetric with respect to $B_{2}$ (see also Figure \ref{symfle}).

\begin{remark}
\label{generalremonmax}
\emph{Given a framework $(G,p)\in \mathscr{R}_{(G,S)}$, we need to specify a type $\Phi:S\to \textrm{Aut}(G)$ in order to apply the symmetry-extended version of Maxwell's rule (Theorem \ref{maxwellsrulewithchar}) to $(G,p)$ and $S$, because $\Phi$ determines the characters $X_{e}$ and $X_{i}$. Of course, we also need to make sure that $\Phi$ is a homomorphism, for otherwise the external and internal representation (with respect to $G$ and $\Phi$) are not matrix representations of $S$ (see Remark \ref{homrepre}).\\\indent The conditions under which $(G,p)\in \mathscr{R}_{(G,S)}$ is of a unique type are given in \cite{BS1}. In the same paper, it is also shown that if $(G,p)\in \mathscr{R}_{(G,S)}$ is of a unique type $\Phi$, then $\Phi$ is a homomorphism.
In particular, this is the case if $(G,p)$ is an injective realization, so that the external and internal representation are uniquely determined in this case and Theorem \ref{maxwellsrulewithchar} can be applied to $(G,p)$ and $S$ in a unique way. Moreover, for injective realizations in $\mathscr{R}_{(G,S)}$, the characters $X_{e}$ and $X_{i}$ can be found in a particularly easy way (without determining the type $\Phi$) by simply examining the geometric positions of the joints and bars of $(G,p)$ (see \cite{cfgsw, BS1, BS4}).\\\indent Since in \cite{FGsymmax} only injective realizations in $\mathbb{R}^{2}$ and $\mathbb{R}^{3}$ are considered, Theorem \ref{maxwellsrulewithchar} includes the symmetrized version of Maxwell's rule given in \cite{FGsymmax} as a special case.\\\indent
If $(G,p)\in \mathscr{R}_{(G,S)}$ is a non-injective realization, then there may exist more than just one homomorphism $\Phi$ for which $(G,p)\in \mathscr{R}_{(G,S,\Phi)}$, in which case we can apply Theorem \ref{maxwellsrulewithchar} to $(G,p)$ and $S$ by using any one of these homomorphisms. It is also possible that there does not exist any homomorphism $\Phi$ so that $(G,p)\in \mathscr{R}_{(G,S,\Phi)}$, in which case we cannot apply the symmetry-extended version of Maxwell's rule to $(G,p)$ and $S$ at all. See again \cite{BS1} for details.}
\end{remark}

\section{Further work}

There exist a number of other classical counting rules, similar to Maxwell's rule, for which symmetry extensions have been derived using techniques from group representation theory (see \cite{FGsymmax, FG3, gsw}, for example). Like the symmetry-extended version of Maxwell's rule in \cite{FGsymmax}, however, these rules are also not presented with a mathematically precise formulation nor with a thorough mathematical foundation or proof. One approach to establish mathematical proofs for these rules is to appropriately modify or extend the methods and results presented in this paper. An alternate approach is presented by J.C. Owen and S.C. Power in \cite{owen}.\\\indent
In \cite{FGsymmax},  a symmetry-extended version of Maxwell's rule for isostatic \emph{pinned} bar and joint frameworks is given. If we define an external and internal representation in the same way as in Section \ref{subsec:intext} (by taking into account only the \emph{unpinned} joints in the definition of the external representation), then we can easily establish a mathematical proof of this rule by slightly modifying the results of this paper. In fact, since a pinned framework is firmly anchored to the ground and hence does not possess any infinitesimal rigid motions, a proof of this rule requires significantly less work than the proof of the symmetry-extended version of Maxwell's rule given in this paper.\\\indent
In \cite{owen}, J.C. Owen and S.C. Power use their analysis of symmetric pinned frameworks as a starting point to establish symmetry-extended counting rules for a variety of other geometric constraint systems.\\\indent
In \cite{gsw}, a symmetrized counting rule for body-bar frameworks is used  to show that an isostatic body-bar framework must satisfy some easily stated restrictions on the number of bodies and bars that are `fixed' by various symmetry operations of the structure. The key to proving these results is to suitably adapt the definitions of the external and internal representation given in Section \ref{subsec:intext}, and then to establish a result analogous to Lemma \ref{replemma}. Once one has shown, based on this lemma, that the rigidity matrix of a body-bar framework can be put into a block-diagonalized form, one needs to follow the steps of Section 4 of this paper, appropriately adapting the results in each step.
\\\indent Proving the symmetry-extended mobility criterion for body-hinge frameworks given in \cite{FG3} is somewhat more complicated. In particular, one has to explain in detail how to associate appropriate group representations to the various hinge constraints.

\section*{Acknowledgements}

We would like to thank Walter Whiteley for his invaluable advice on shaping this presentation, as well as Simon Guest for fruitful and interesting discussions at the BIRS workshop in July 2008.

\end{document}